\documentclass[a4paper,11pt]{amsart}

\usepackage{graphicx}
\usepackage{mathrsfs} % required for \mathscr font

\def\ie{i\textup.e\textup.}

\def\rank{\operatorname{rank}}

\newtheorem*{named_theorem}{Theorem}
\newtheorem*{theorem}{Theorem}
\newtheorem*{main_lemma}{Main Lemma}
\newtheorem{property}{Property}
\newtheorem{proposition}{Proposition}[section]
\newtheorem{lemma}{Lemma}[section]
\newtheorem{corollary.proposition}{Corollary}[proposition]

\theoremstyle{definition}
\newtheorem*{definition}{Definition}

\theoremstyle{remark}
\newtheorem{remark}{Remark}[section]

\numberwithin{equation}{section}

\begin{document}

%% AMS-LaTeX
%% Paper 2 - Version GAMMA
%% Top Matter
%% (Title and so on...)
%%%%%%%%%%%%%%%%%%%%%%%%%%%%%%%%%%%%%%%%%%%%%%%%%%%%%%%%%%%%%%%%%%%%%%%%%%%

%!TEX root = paper2-2005.tex

\title[On torsion-free groups with finite regular file bases]
{On torsion-free groups with\\ finite regular file bases}

\author{Alexey Muranov}

\address{Dept of Mathematics, Vanderbilt University,
1326 Stevenson Center,
Nashville, TN 37240--0001, USA.}

\curraddr{Institut Camille Jordan\\
Universit\'e Lyon 1\\
43 blvd du 11 novembre 1918\\
69622 Villeurbanne cedex\\
France.}

\email{muranov@math.univ-lyon1.fr}

\thanks{This work was supported in part by the NSF grant DMS 0245600
of Alexander Ol'shanskiy and Mark Sapir.}



\date{\today}

%\dedicatory{}

\commby{Alexander Kleshchev}

\begin{abstract}
The following question was asked by V.~V.~Bludov in 
The Kourov\-ka Notebook in 1995:
If a torsion-free group $G$ has a finite system of generators 
$a_1$, \dots, $a_n$
such that every element of $G$ has a unique presentation in the form
$a_1^{k_1}\dots a_n^{k_n}$ where $k_i\in\mathbb Z$,
is it true that $G$ is virtually polycyclic?
The answer is ``not always.''
A counterexample is constructed in this paper as a group presented
by generators and defining relations.
\end{abstract}

\maketitle

\tableofcontents

%%%%%%%%%%%%%%%%%%%%%%%%%%%%%%%%%%%%%%%%%%%%%%%%%%%%%%%%%%%%%%%%%%%%%%%%%%%
 %<------------\|- Top matter
%% AMS-LaTeX
%% Paper 2 - Version GAMMA
%% Part 01
%% (Introduction)
%%%%%%%%%%%%%%%%%%%%%%%%%%%%%%%%%%%%%%%%%%%%%%%%%%%%%%%%%%%%%%%%%%%%%%%%%%%

%!TEX root = paper2-2005.tex

\section{Introduction}
\label{section:introduction}

\begin{definition}
A group is \emph{polycyclic\/} if
%and only if
it has a finite subnormal series with cyclic factors.
A group is \emph{virtually polycyclic\/} if
%and only if
it has a polycyclic subgroup of finite index.
\end{definition}

Subgroups and homomorphic images of polycyclic and virtually polycyclic 
groups are polycyclic or, respectively, virtually polycyclic.
Finitely generated abelian groups are polycyclic.
Free groups of rank at least $2$ are not virtually polycyclic.

Vasiliy Bludov posed the following question in 1995
(Problem~13.11 in \cite{kn13:1995:knupgt-eng}; 
see also Question~6 in~\cite{Bludov:1995:fbg-eng}):
\begin{quote}
If a torsion-free group $G$ has a finite system of generators 
$a_1$, \dots, $a_n$ such that every element of $G$ has a unique 
presentation in the form $a_1^{k_1}\dots a_n^{k_n}$ 
where $k_i\in\mathbb Z$,
is it true then that $G$ is virtually polycyclic?
\end{quote}
The following theorem answers Bludov's question negatively:

\begin{theorem}
Provided\/ $n\ge63$\textup,
there exists a group\/ $G$ and pairwise
distinct elements\/ $a_1$\textup, \dots\textup, $a_n$ in\/ $G$
such that\/\textup:
\begin{enumerate}
\setcounter{enumi}{-1}
\item
	$G$ is generated by\/ $\{a_1,\dots,a_n\}$\textup;
\item
	every\/ $n-21$ elements out of\/ $\{a_1,\dots,a_n\}$ 
	freely generate a free subgroup such that 
	every two elements of this subgroup\/ $F$ 
	are conjugate in\/ $G$ only if they are conjugate in\/ $F$ itself\/
	\textup(in particular\textup, 
	$G$ is not virtually polycyclic\/\textup{);}
\item
\label{item:file_basis_in_main_theorem}
	for every element\/ $g$ of\/ $G$\textup,
	there is a unique\/ $n$-tuple\/ $(k_1,\dots,k_n)\in\mathbb Z^n$ 
	such that\/ $g=a_1^{k_1}\dots a_n^{k_n}$\textup;
\item
	$G$ is torsion-free\/\textup;
\item
\label{item:direct_limit_in_main_theorem}
	$G$ is the direct limit of a sequence
	of hyperbolic groups with respect to a family of surjective 
	homomorphisms\/\textup;
\item
	$G$ is recursively presented and has solvable word and conjugacy 
	problems\textup.
\end{enumerate}
\end{theorem}

(The \emph{direct limit\/} in condition 
(\ref{item:direct_limit_in_main_theorem}) 
is the direct limit in the sense of Bourbaki; see ``Algebra~I.'')

This theorem is proved in the end of 
Section~\ref{section:properties__group} of this paper.

In terms of \emph{file bases\/} 
(see definitions in \cite{Bludov:1995:fbg-eng}),
item (\ref{item:file_basis_in_main_theorem}) 
of the theorem states that the elements $a_1$, \dots, $a_n$ 
in this order form a \emph{regular\/} file basis of~$G$.
The argument below shows that every group with a $2$-element file basis 
is polycyclic
(therefore, the assumption ``$n\ge63$'' cannot be replaced with 
``$n\ge2$'').

Every group that has an $n$-element file basis is 
\emph{$n$-boundedly generated\/}
(\ie, is the product of an $n$-tuple of its cyclic subgroups).
The following two theorems 
(see \cite{Ito:1955:UdPvzaG-ger} and \cite{LennoxRoseblade:1980:sppg}) 
together yield that every $2$-boundedly generated group is polycyclic:

\begin{named_theorem}[It\^o, 1955]
If a group\/ $G$ has abelian subgroups\/ $A$ and\/ $B$ such that\/ 
$G=AB$\textup, 
then\/ $G$ is metabelian\textup.
\end{named_theorem}

\begin{named_theorem}[Lennox, Roseblade, 1980]
If a soluble group\/ $G$ has polycyclic subgroups\/ $A$ and\/ $B$ 
such that\/ $G=AB$\textup, 
then\/ $G$ is polycyclic\textup.
\end{named_theorem}

Some examples of non-polycyclic groups
with finite regular file bases and with torsion are known.
It is shown in \cite{Bludov:1995:fbg-eng} (see Example~3 therein) that
the semidirect product of the rank-$n$ free abelian group $\mathbb Z^n$
with the symmetric group $S_n$ defined via the action of $S_n$ on 
$\mathbb Z^n$ by naturally permuting the standard generators
has an $n$-element regular file basis of infinite-order elements.
If $n\ge5$, then such a semidirect product is virtually polycyclic but 
not polycyclic.

The primary tool used to establish properties of the group 
(or rather the series of groups) constructed in 
this paper is the Main Theorem in~\cite{Muranov:2005:dsmcbgbsg}.
\par

%%%%%%%%%%%%%%%%%%%%%%%%%%%%%%%%%%%%%%%%%%%%%%%%%%%%%%%%%%%%%%%%%%%%%%%%%%%
 %<------------\|- Introduction
%% AMS-LaTeX
%% Paper 2 - Version GAMMA
%% Part 02
%% (Group Construction)
%%%%%%%%%%%%%%%%%%%%%%%%%%%%%%%%%%%%%%%%%%%%%%%%%%%%%%%%%%%%%%%%%%%%%%%%%%%

%!TEX root = paper2-2005.tex

\section{Group construction}
\label{section:group_construction}

Take an arbitrary integer $n\ge63$.
Choose a positive $\lambda_1<1$ such that
\begin{equation}
\label{display:(2.1)}
\Bigl(4+\frac{2n\lambda_1}{1-\lambda_1}\Bigr)\lambda_1\le\frac{1}{n}.
\end{equation}
(It suffices if, for example, $0<\lambda_1\le1/(5n)$.)

Let $\mathfrak A=\{x_1,\dots,x_n\}$ be an $n$-letter alphabet.
Call a group word $w$ over $\mathfrak A$ \emph{regular\/} if
%and only if
it has the form $x_1^{k_1}x_2^{k_2}\dots x_n^{k_n}$ where 
$k_i\in\mathbb Z$.
(In particular, regular group words are cyclically reduced.)
Call a group word $w$ \emph{counter-regular\/} if
%and only if
$w^{-1}$ is regular.
The empty word is both regular and counter-regular.
So is every group word that is a \emph{letter power}.

Impose an order on the set of all reduced group words in the alphabet 
$\mathfrak A$ to make it order-isomorphic with the set of natural numbers.
For example, use the \emph{deg-lex\/} order: 
$x_1$, $x_1^{-1}$, $x_2$, \dots, 
$x_n^{-1}$, $x_1^2$, $x_1x_2$, $x_1x_2^{-1}$, \dots.
Define a sequence of sets $\{\mathcal R_i\}_{i=0}^{+\infty}$ inductively.

First, let $\mathcal R_0=\emptyset$.

Second, if $i>0$ and every group word over $\mathfrak A$ equals 
a regular word
modulo the relations $r=1$, $r\in\mathcal R_{i-1}$, then
let $\mathcal R_i=\mathcal R_{i-1}$.

Last, if $i>0$ and some group word over $\mathfrak A$ is not equal to any 
regular word modulo the relations $r=1$, $r\in\mathcal R_{i-1}$, then
let $w_i$ be the least reduced group word 
(with respect to the chosen order)
that does not start with $x_1^{\pm1}$, does not end with $x_n^{\pm1}$,
and is not equal to any regular word modulo the relations 
$r=1$, $r\in\mathcal R_{i-1}$ (clearly such $w_i$ exists).
Choose a natural number $m_i$ such that
\begin{gather}
%\spreadlines{2\jot}
\label{display:(2.2)}
i>j\,\Rightarrow\,m_i\ne m_j,\\[2\jot]
\label{display:(2.3)}
i>j\,\Rightarrow\,nm_i+|w_i|\ge nm_j+|w_j|,\\[2\jot]
\label{display:(2.4)}
\lambda_1\bigl(nm_i+|w_i|\bigr)\ge|w_i|.
\end{gather}
(These conditions may be satisfied by choosing a sufficiently large~$m_i$.)
Define
\begin{equation}
\label{display:(2.5)}
r_i=x_1^{m_i}x_2^{m_i}\dots x_n^{m_i}w_i^{-1},
\end{equation}
and let $\mathcal R_i=\mathcal R_{i-1}\cup\{r_i\}$.
Note that $|r_i|=nm_i+|w_i|$.

Eventually, let
\begin{equation}
\label{display:(2.6)}
\mathcal R=\bigcup_{i=1}^{+\infty}\mathcal R_i,
\end{equation}
and let $G$ be the group defined by the presentation
$\langle\,\mathfrak A\,\|\,r=1,r\in\mathcal R\,\rangle$.

Observe that all elements of $\mathcal R$ are cyclically reduced.

Along with $\{\mathcal R_i\}_{i=0}^{+\infty}$, the sequences
$\{w_i\}_{i=1,\dots}$, $\{m_i\}_{i=1,\dots}$, and $\{r_i\}_{i=1,\dots}$,
infinite or finite, have been constructed.

Note that the above construction scheme does not specify $\mathcal R$ 
uniquely and is flexible as to the choice of 
$\{w_i\}_{i=1,\dots}$ and $\{m_i\}_{i=1,\dots}$.
\par

%%%%%%%%%%%%%%%%%%%%%%%%%%%%%%%%%%%%%%%%%%%%%%%%%%%%%%%%%%%%%%%%%%%%%%%%%%%
 %<------------\|- Group construction
%% AMS-LaTeX
%% Paper 2 - Version GAMMA
%% Part 03
%% (Definitions)
%%%%%%%%%%%%%%%%%%%%%%%%%%%%%%%%%%%%%%%%%%%%%%%%%%%%%%%%%%%%%%%%%%%%%%%%%%%

%!TEX root = paper2-2005.tex

\section{Definitions}
\label{section:definitions}

Most definitions needed for reading this paper may be found 
in~\cite{Muranov:2005:dsmcbgbsg}.
Standard terms are also defined 
in \cite{LyndonSchupp:2001:cgt}, \cite{ChiColHue:1981:agp}, 
\cite{Olshanskii:1991:gdrg-eng} in similar ways.
For instance, geometric definitions of \emph{maps\/} and \emph{diagrams}
are given in \cite{Olshanskii:1991:gdrg-eng}, and combinatorial
ones in~\cite{Muranov:2005:dsmcbgbsg}.
Some definitions are reviewed in this section.

If $X$ is a set, then $\|X\|$ denotes the cardinality of~$X$.

A \emph{graph\/} $\Gamma$ is an ordered pair ${(}\Gamma(0),\Gamma(1){)}$ 
where $\Gamma(0)$ is the set of \emph{vertices\/} of $\Gamma$,
and $\Gamma(1)$ is its set of \emph{edges}.
It will always be assumed that $\Gamma(0)$ is nonempty.
In the context of this paper, \emph{graph\/} is a synonym of 
\emph{combinatorial\/ $1$-complex}.
A \emph{combinatorial\/ $2$-complex\/} $\Phi$ is an ordered triple
${(}\Phi(0),\Phi(1),\Phi(2){)}$ where $\Phi(0)$, $\Phi(1)$, and $\Phi(2)$
are its sets of vertices, edges, and \emph{faces\/}, respectively.

%The term \emph{graph\/} in this paper refers to an object that in graph 
%theory is called an \emph{undirected pseudograph\/} 
%(multiple edges and loops are admissible).
%If $\Gamma$ is a graph, then $\Gamma(0)$ denotes its (nonempty) set of 
%\emph{vertices},
%and $\Gamma(1)$ denotes its set of \emph{edges}.
%The graph itself is an ordered pair:
%$\Gamma={(}\Gamma(0),\Gamma(1){)}$.
%In the context of this paper, \emph{graph\/} is a synonym of 
%\emph{combinatorial\/ $1$-complex}.
%\emph{Combinatorial\/ $2$-complexes\/} also play an important part in 
%this paper.
%If $\Phi$ is a $2$-complex, then $\Phi(0)$, $\Phi(1)$, and $\Phi(2)$
%denote its sets of vertices, edges, and \emph{faces\/} (cells), 
%respectively.
%The $2$-complex $\Phi$ itself is the ordered triple
%${(}\Phi(0),\Phi(1),\Phi(2){)}$.

If $v$ is a vertex of a complex $\Gamma$, then
the number of all edges of $\Gamma$ incident to $v$ plus the number
of all \emph{loops\/} incident to $v$ (thus loops are counted twice)
is called the \emph{degree\/} of the vertex $v$ 
and is denoted by $d(v)$ or $d_\Gamma(v)$.

Here follow several definitions related to \emph{paths\/}
in combinatorial complexes.

%A \emph{path\/} is a finite sequence of alternating vertices and
%oriented edges such that the following conditions hold:
%it starts with a vertex and ends with a vertex;
%the vertex immediately preceding an oriented edge is its \emph{tail\/};
%the vertex succeeding an oriented edge is its \emph{head}.
%%The \emph{initial vertex\/} of a path is its first vertex,
%%the \emph{terminal vertex\/} of a path is its last vertex;
%%together they are the \emph{end-vertices\/} of the path.
%The \emph{length\/} of a path
%$p={(} v_0,e_1,v_1,\dots,e_n,v_n{)}$ is $n$;
%it is denoted by~$|p|$.
%The vertices $v_1,\dots,v_{n-1}$ of this path are called
%its \emph{intermediate\/} vertices.
%A \emph{trivial path\/} is a path of length zero.
%By abuse of notation, a path of the form ${(} v_1,e,v_2{)}$,
%where $e$ is an oriented edge from $v_1$ to $v_2$, shall be denoted 
%by $e$, and a trivial path ${(} v{)}$ shall be denoted by~$v$.
%
%The \emph{inverse path\/} to a path $p$ is defined naturally and
%is denoted by~$p^{-1}$.
%If the terminal vertex of a path $p_1$ coincides with
%the initial vertex of a path $p_2$, then the \emph{product\/}
%$p_1 p_2$ is defined.
%%A path $s$ is called an \emph{initial subpath\/} of a path $p$ if
%%%and only if
%%$p=sq$ for some path~$q$.
%%A path $s$ is called a \emph{terminal subpath\/} of a path $p$ if
%%%and only if
%%$p=qs$ for some path~$q$.

The length of a path $p$ is denoted by~$|p|$.
A \emph{cyclic path\/} is a path whose \emph{terminal\/} and \emph{initial\/}
vertices coincide.
A \emph{cycle\/} is the set of all \emph{cyclic shifts\/} of a cyclic path.
The cycle represented by a cyclic path $p$ shall be denoted 
by~$\langle p \rangle$.
The \emph{length\/} of a cycle $c$, denoted by $|c|$,
is the length of an arbitrary representative of~$c$.
A \emph{trivial cycle\/} is a cycle of length zero.
A path $p$ is a \emph{subpath of a cycle\/} $c$ if
%and only if
for some representative $r$ of $c$ and for some natural $n$,
$p$ is a subpath of $r^n$ (\ie, of the product of $n$ copies of~$r$).

A path is \emph{reduced\/} if
%and only if
it does not have a
subpath of the form $ee^{-1}$ where $e$ is an oriented edge.
A cyclic path is \emph{cyclically reduced\/} if
%and only if
it is reduced, and its first
oriented edge is not inverse to its last oriented edge.
(For example, all trivial paths are cyclically reduced.)
A cycle is \emph{reduced\/} if
%and only if
it consists of cyclically reduced cyclic paths.
A path is \emph{simple\/} if
%and only if
it is nontrivial, reduced, and
none of its \emph{intermediate\/} vertices appears in it more than once.
A cycle is \emph{simple\/} if
%and only if
it consists of simple cyclic paths.
A path is \emph{maximal\/} in some set of paths if
%and only if
it is not a proper subpath of any other path in this set.

An \emph{oriented arc\/} in a complex $\Gamma$ is a simple path
whose all intermediate vertices have degree $2$ in~$\Gamma$.
A \emph{non-oriented arc}, or simply an \emph{arc}, is defined as a pair
of mutually inverse oriented arcs.
%The \emph{length\/} of a non-oriented arc $u$ is the length
%of either of the two associated oriented arcs and is denoted by~$|u|$.
%The concepts of a \emph{subarc}, a \emph{maximal arc},
%and similar, are self-explanatory.
Edges and oriented edges shall be viewed as special cases of 
arcs and oriented arcs, respectively.
%The \emph{set of edges\/} of an oriented arc (respectively of an arc) 
%is the set of all
%edges that with some orientation occur in that oriented arc
%(respectively in either of the two associated oriented arcs).
%An \emph{intermediate vertex\/} of an arc is an intermediate vertex
%of either of the associated oriented arcs.
Two arcs, or oriented arcs, are said to \emph{overlap\/} if
%and only if
they have common edges.
A set of arcs, or oriented arcs, is called \emph{non-overlapping\/} if
%and only if
no two distinct elements of this set overlap.
An arc $u$ \emph{lies\/} on a path $p$ if
%and only if
at least one of the oriented arcs associated with $u$ is a subpath of~$p$.

The remaining definitions concern maps and diagrams.

Every oriented face of a $2$-complex has a unique \emph{boundary cycle}.
Every (non-oriented) face of an oriented combinatorial surface has 
a unique \emph{contour},
which is defined as the boundary cycle of the corresponding oriented face.
The contour of a face $\Pi$ shall be denoted by~$\partial\Pi$.
Thus the formula $\langle p \rangle=\partial\Pi$ means that the path $p$
is a representative of the contour of the face~$\Pi$.
The \emph{degree\/} of a face $\Pi$ is the length of the contour of $\Pi$, 
denoted by $|\partial\Pi|$.

A \emph{map\/} is a connected subcomplex
of an oriented combinatorial sphere
together with a set of cycles called its \emph{contours\/}
(see \cite{Muranov:2005:dsmcbgbsg} for a precise definition).
In every map, each oriented edge appears either in the contour of some 
face of the map or in a contour of the map.
The union of all contours of a map $\Delta$ shall be denoted 
by~$\partial\Delta$.
A map without faces is called \emph{degenerate}.
A map is \emph{simple\/} if
%and only if
all of its contours are simple cycles without common vertices.
A \emph{disc map\/} is a map with exactly one contour.
A disc map is simple if and only if
its underlying $2$-complex is a \emph{combinatorial disc}.
An \emph{annular map\/} is a map with exactly two contours.

The contours of faces in a disc map shall be thought of as
\emph{oriented counterclockwise}, and the contour of the map itself
shall be thought of as \emph{oriented clockwise\/}
(this is consistent with the way disc maps are pictured in this paper).

A \emph{selection\/} on a face $\Pi$ of a map $\Delta$ is a set of 
nontrivial reduced subpaths of the contour of $\Pi$ closed under taking 
nontrivial subpaths.
A \emph{selection\/} on a map $\Delta$ is a set of nontrivial reduced 
subpaths of the contours of faces of $\Delta$ closed under taking 
nontrivial subpaths.
An \emph{S-map} is a map with selection.
A path in an S-map is \emph{selected\/} if
%and only if
it belongs to the selection.
A path $p$ is \emph{double-selected\/} if
%and only if
both $p$ and $p^{-1}$ are selected.
An arc is \emph{double-selected\/} if
%and only if
both corresponding oriented arcs are selected.
An external arc is \emph{selected\/} if
%and only if
one of the corresponding oriented arcs is selected.

A \emph{diagram\/} over a presentation 
$\langle\,\mathfrak B\,\|\,\mathcal S\,\rangle$
is a map together with a function which labels every oriented edge of this
map with a letter from $\mathfrak B^{\pm1}$ so that mutually inverse 
oriented edges are labelled
with mutually inverse group letters, and the contour of every face 
has a representative ``labelled'' with an element of~$\mathcal S^{\pm1}$.
A diagram together with a selection on the underlying map is called 
an \emph{S-diagram}.

In a diagram over $\langle\,\mathfrak B\,\|\,\mathcal S\,\rangle$,
let the \emph{label\/} of an oriented edge be the element of 
$\mathfrak B^{\pm1}$
that labels this oriented edge, and let the \emph{label\/} of a path
be the group word that ``reads'' on this path.
Let the \emph{label\/} of a (non-oriented) edge be the letter 
$x$ of $\mathfrak B$ such that $x^{+1}$ and $x^{-1}$ label 
the two corresponding oriented edges.
If $q$ is an oriented edge or a path, the label $q$ shall be denoted 
by $\ell (q)$.

Sometimes it is convenient to distinguish between a letter $x$ of 
the alphabet $\mathfrak B$
and the corresponding letter of the group alphabet $\mathfrak B^{\pm1}$ 
denoted by $x^{+1}$ or simply by~$x$.
Call letters of the alphabet $\mathfrak B$ \emph{basic letters\/}
and letters of the group alphabet $\mathfrak B^{\pm1}$ 
\emph{group letters}.
Thus, the basic letter $x_1$, the group letter $x_1$,
and the one-letter group word $x_1$ are three different things.
The labels of edges of a diagram are basic letters,
but the labels of oriented edges are group letters.
Say that a group word $w$ has a basic letter $x$ in it if
%and only if
at least one of the two group letters $x^{\pm1}$ occurs in~$w$.

The \emph{mirror copy\/} of a map or diagram is the map or diagram, 
respectively, with the same underlying $2$-complex, and the same 
labelling function in the case of diagram,
but with the opposite orientation of all the faces and contours.

According to van Kampen's lemma (see \cite{Olshanskii:1991:gdrg-eng}),
a relation $w=1$ is a consequence of a set of relations 
$\{\,r=1\mid r\in\mathcal S\,\}$
if and only if there exists a disc diagram $\Delta$ over 
$\langle\,\mathfrak B\,\|\,\mathcal S\,\rangle$, where $\mathfrak B$ is 
an appropriate alphabet,
such that the label of some representative of $\partial\Delta$ is~$w$.
Such a diagram is called a \emph{deduction diagram\/} for 
the relation $w=1$.

Call a pair of distinct faces $\{\Pi_1,\Pi_2\}$ in a diagram $\Delta$
\emph{immediately cancellable\/} if
%and only if
there exist paths $p_1$ and $p_2$ with a nontrivial common initial subpath
such that $\langle p_1\rangle=\partial\Pi_1$,
$\langle p_2^{-1}\rangle=\partial\Pi_2$, and $\ell (p_1)=\ell (p_2)$.
A diagram without immediately cancellable pairs of faces is called 
\emph{weakly reduced}.
(In \cite{Muranov:2005:dsmcbgbsg} immediately cancellable pairs of faces 
and weakly reduced 
diagrams are called simply \emph{cancellable\/} and \emph{reduced}, 
respectively.)
If there exists a deduction diagram for a relation $w=1$ over 
a presentation
$\langle\,\mathfrak B\,\|\,\mathcal S\,\rangle$,
then such a diagram may be chosen weakly reduced,
and even reduced in the sense of \cite{Olshanskii:1991:gdrg-eng},
which is a stronger property
(take an arbitrary deduction diagram and ``eliminate'' immediately 
cancellable pairs).

Call a group presentation $\langle\,\mathfrak B\,\|\,\mathcal S\,\rangle$
\emph{strongly aspherical\/} if
%and only if
no spherical diagram over $\langle\,\mathfrak B\,\|\,\mathcal S\,\rangle$
is weakly reduced.
\par

%%%%%%%%%%%%%%%%%%%%%%%%%%%%%%%%%%%%%%%%%%%%%%%%%%%%%%%%%%%%%%%%%%%%%%%%%%%
 %<------------\|- Definitions
%% AMS-LaTeX
%% Paper 2 - Version GAMMA
%% Part 04
%% (The Main Lemma)
%%%%%%%%%%%%%%%%%%%%%%%%%%%%%%%%%%%%%%%%%%%%%%%%%%%%%%%%%%%%%%%%%%%%%%%%%%%

%!TEX root = paper2-2005.tex

\section{The Main Lemma}
\label{section:the_main_lemma}

This paper relies on a technical result of \cite{Muranov:2005:dsmcbgbsg} 
called there ``the Main Theorem.''
In fact, it suffices to use a simplified form of that result, 
which is stated in this section as the Main Lemma.

``Inductive Lemmas'' and ``the Main Theorem''
in this section refer to the corresponding statements 
in~\cite{Muranov:2005:dsmcbgbsg}.
Because they are not used outside of this section, because their 
statements 
involve properties of maps that are not defined in this paper,
and because the Main Lemma may be viewed as a particular case of 
the Main Theorem, the Main Theorem and Inductive Lemmas are not 
stated here.

\begin{definition}
Call a map $\Delta$ \emph{semisimple\/} if
%and only if
every edge of $\Delta$ is incident to a face of~$\Delta$.
\end{definition}

\begin{definition}
Let $\Delta$ be a semisimple S-map.
Let $S$ be the number of selected external edges of~$\Delta$.
Let $\mu$ be a real number.
Then the S-map $\Delta$ is said to satisfy the \emph{condition\/} 
$\mathcal X(\mu)$ if
%and only if
$$
S\ge\|\Delta(1)\|-\mu\!\sum_{\Pi\in\Delta(2)}\!|\partial\Pi|.
$$
\end{definition}

Note that the definition of the condition $\mathcal X$ in 
\cite{Muranov:2005:dsmcbgbsg}, 
given there only for simple S-maps, 
is extended here to all semisimple S-maps.

\begin{proposition}
\label{proposition:(4.1)}
Let\/ $\Delta$ be an S-map\textup, and let\/ $\mu$ be a real number\textup.
Suppose every maximal semisimple S-submap of\/ $\Delta$
satisfies the condition\/ $\mathcal X(\mu)$\textup.
Let\/ $S$ be the number of selected external edges of\/ $\Delta$\textup.
Then
$$
S\ge(1-2\mu)\!\sum_{\Pi\in\Delta(2)}\!|\partial\Pi|.
$$
\end{proposition}

\begin{proof}
Consider an arbitrary maximal semisimple S-submap $\Delta'$ of~$\Delta$.
Let $S'$ be the number of selected external edges of~$\Delta'$.
Let $\Sigma'$ be the sum of the degrees of all the faces of~$\Delta'$.
Then
$$
S'=2S'-S'\ge2\|\Delta'(1)\|-2\mu\Sigma'-S'
\ge S'+\Sigma'-2\mu\Sigma'-S'=(1-2\mu)\Sigma'.
$$

Now, consider $\Delta$ itself.
Remove from the underlying $2$-complex of $\Delta$
all edges that are not incident to any face.
All connected components of the remaining $2$-complex
(some, possibly, consisting of just one vertex)
define semisimple S-submaps of~$\Delta$.
These S-submaps are precisely all maximal semisimple S-submaps of~$\Delta$.
Therefore, every selected external edge of $\Delta$, 
as well as every face of $\Delta$,
belongs to exactly one maximal semisimple S-submap of~$\Delta$.
Hence, the desired inequality for $\Delta$ may be obtained by
adding up the analogous inequalities for all of its maximal 
semisimple S-submaps.
\end{proof}

Although it takes little effort to derive the Main Lemma formulated below 
from the Main Theorem, it will be even easier if a modified version 
of the Main Theorem with slightly weakened hypotheses is used instead.
It can be observed that the Main Theorem will still hold true if 
``simple S-map'' is replaced with ``semisimple S-map'' in the hypotheses.
Indeed, analysis of the proofs of Inductive Lemmas shows that
their statements can be modified in the following way 
without turning them false:
\begin{enumerate}
\item
	The hypotheses of Inductive Lemma~1 may be weakened by requiring 
	the condition $\mathcal X(\mu)$ only from all simple disc 
	S-submaps of $\Delta$ all of whose edges are internal in $\Delta$,
	rather than from all proper simple disc S-submaps.
\item
	The hypotheses of Inductive Lemma~2 may be weakened by requiring 
	the S-map $\Delta$ to be only semisimple, rather than simple.
\end{enumerate}
The original proofs with a few obvious changes apply to the modified lemmas.

The statement obtained from the Main Theorem by replacing
``simple'' with ``semisimple'' shall be referred to as 
the Modified Main Theorem.
The original proof of the Main Theorem applies to the Modified Main Theorem
if Inductive Lemmas, used in the proof, are modified as described.
(The condition $\mathcal X$ in the conclusions of the Modified Main Theorem 
and modified Inductive Lemma~2 shall be understood in the above-defined 
broader sense.)

\begin{remark}
The Modified Main Theorem can also be obtained as a corollary of 
the Main Theorem by ``cutting'' a given semisimple S-map into its 
``simple components.''
\end{remark}

\begin{definition}
Let $\Delta$ be an S-map.
Let $\lambda_1,\lambda_2\in[0,1]$.
The S-map $\Delta$ is said to satisfy the \emph{condition\/}
$\mathcal B(\lambda_1,\lambda_2)$ if
%and only if
it satisfies the following three conditions:
\begin{description}
\item[$\mathcal B_0$]
	For each face $\Pi$ of $\Delta$, the contour of $\Pi$ has
	at least one selected subpath and at most one maximal selected subpath
	(note that if all nontrivial subpaths are selected,
	then there is no maximal selected subpath).

\item[$\mathcal B_1(\lambda_1)$]
	For each face $\Pi$ of $\Delta$,
	there is a selected subpath of $\partial\Pi$
	of length at least $(1-\lambda_1)|\partial\Pi|$.

\item[$\mathcal B_2(\lambda_2)$]
	For each face $\Pi$ of $\Delta$,
	the length of every double-selected arc incident to $\Pi$
	is at most $\lambda_2|\partial\Pi|$.
\end{description}
A map $\Delta$ is said to satisfy the condition 
$\mathcal B(\lambda_1,\lambda_2)$ if
%and only if
there exists a selection on $\Delta$ such that $\Delta$ with this selection
satisfies $\mathcal B(\lambda_1,\lambda_2)$.
\end{definition}

For all admissible values of $\lambda_1$ and $\lambda_2$,
the condition $\mathcal B(\lambda_1,\lambda_2)$ is equivalent to
the condition $\mathcal A(1,1,1;\lambda_1,\lambda_2,0,0)$ 
in~\cite{Muranov:2005:dsmcbgbsg}.

\begin{main_lemma}
Let\/ $\Delta$ be an S-map with at most\/ $3$ contours\textup.
Let\/ $\lambda_1,\lambda_2\in[0,1]$\textup.
Suppose\/ $2\lambda_1+13\lambda_2<1$
and\/ $\Delta$ satisfies\/ $\mathcal B(\lambda_1,\lambda_2)$\textup.
Let\/ $S$ be the number of selected external edge of\/ $\Delta$\textup,
and let\/ $\mu=\lambda_1+5\lambda_2$\textup.
Then
$$
S\ge(1-2\mu)\!\sum_{\Pi\in\Delta(2)}\!|\partial\Pi|.
$$
\end{main_lemma}

\begin{proof}
Every maximal semisimple S-submap of $\Delta$ has at most $3$ contours
(see the description of maximal semisimple S-submaps given in the proof 
of Proposition~\ref{proposition:(4.1)})
and satisfies the condition $\mathcal B(\lambda_1,\lambda_2)$.
By the Modified Main Theorem, and Remarks 6.3 and 6.4 
in \cite{Muranov:2005:dsmcbgbsg},
all maximal semisimple S-submaps of $\Delta$ satisfy 
the condition $\mathcal X(\mu)$.
The desired inequality now follows from 
Proposition~\ref{proposition:(4.1)}.
\end{proof}

%%%%%%%%%%%%%%%%%%%%%%%%%%%%%%%%%%%%%%%%%%%%%%%%%%%%%%%%%%%%%%%%%%%%%%%%%%%
 %<------------\|- The Main Lemma
%% AMS-LaTeX
%% Paper 2 - Version GAMMA
%% Part 05
%% (Properties of the group)
%%%%%%%%%%%%%%%%%%%%%%%%%%%%%%%%%%%%%%%%%%%%%%%%%%%%%%%%%%%%%%%%%%%%%%%%%%%

%!TEX root = paper2-2005.tex

\section{Properties of the group}
\label{section:properties__group}

In this section some properties of the above-defined group $G$
are established.
In particular, it is shown that $G$ is an example that provides
the negative answer to Bludov's question.

Adopt the notation of Section~\ref{section:group_construction}.
In particular, define $n$ and $\lambda_1$ the same way.
Let $\lambda_2=2/n$, $\mu=\lambda_1+5\lambda_2$.
Then
\begin{equation}
\label{display:(5.1)}
2\lambda_1+13\lambda_2<1
\end{equation}
and
\begin{equation}
\label{display:(5.2)}
1-2\mu-2\lambda_1-\frac{2n\lambda_1^2}{1-\lambda_1}\ge1-\frac{21}{n}.
\end{equation}
In particular $\mu<1/2$.
It is well known that $21\times3=63$; this explains why~$63$.

Inequality \thetag{\ref{display:(2.4)}} is equivalent to
\begin{equation}
\label{display:(5.3)}
m_i\ge\frac{1-\lambda_1}{n}\bigl(nm_i+|w_i|\bigr).
\end{equation}
Combining \thetag{\ref{display:(2.4)}}, \thetag{\ref{display:(5.3)}}, and \thetag{\ref{display:(2.3)}}, obtain
\begin{equation}
\label{display:(5.4)}
j\le i\,\Rightarrow\,|w_j|\le\frac{n\lambda_1}{1-\lambda_1}m_i.
\end{equation}
Inequality \thetag{\ref{display:(2.1)}} implies that $\lambda_1<1/(4n)<1/(2n+1)$.
It follows from this and from \thetag{\ref{display:(5.4)}} that
\begin{equation}
\label{display:(5.5)}
j\le i\,\Rightarrow\,|w_j|<\frac{m_i}{2}.
\end{equation}

If $w$ is a group word over $\mathfrak A$,
let $[w]_G$, $[w]_{\mathcal R}$, or simply $[w]$,
denote the element of $G$ represented by~$w$.
Let $a_1$, \dots, $a_n$ be the elements of $G$ represented by
the one-letter group words $x_1$, \dots, $x_n$, respectively
(in terms of ``brackets,'' $a_i=[x_i]$).

In this section a selection on a diagram $\Delta$ is called
\emph{special\/} if
%and only if
the contour of every face $\Pi$ of $\Delta$
has two subpaths $s$ and $t$ such that:
\begin{itemize}
\item
	$\langle st\rangle=\partial\Pi\mathstrut$;
\item
	$s$ is the only maximal selected subpath of
	$\partial\Pi\mathstrut$;
\item
	there is $m\in\mathbb N$ such that either
	$\mathstrut\ell (s)=x_1^{m}x_2^{m}\dots x_n^{m}$ or
	$\ell (s)=x_n^{-m}x_{n-1}^{-m}\dots x_{1}^{-m}\mathstrut$;
\item
	$\displaystyle|s|
	>\frac{\mathstrut n}{\mathstrut 2n-2}|\partial\Pi|$ \ 
	(note that $|\partial\Pi|=|s|+|t|$).
\end{itemize}

On every diagram $\Delta$
over $\langle\,\mathfrak A\,\|\,\mathcal R\,\rangle$,
there exists a unique special selection.
(In fact, if $\Delta$ is a diagram over an arbitrary group presentation,
and a special in the above sense selection on $\Delta$ exists,
then it is unique.)
A diagram over $\langle\,\mathfrak A\,\|\,\mathcal R\,\rangle$
together with a special selection shall be called a
\emph{special\/} S-diagram.
Note that if $s$ is the maximal selected subpath of
the contour of a face $\Pi$ of a special S-diagram
over $\langle\,\mathfrak A\,\|\,\mathcal R\,\rangle$,
then $|s|\ge(1-\lambda_1)|\partial\Pi|$,
which is stronger than the inequality in the definition of
a special selection.

If $\Delta$ is a diagram over
$\langle\,\mathfrak A\,\|\,\mathcal R\,\rangle$,
and $\Pi$ is a face of $\Delta$, then
define the \emph{rank\/} of the face $\Pi$ to be such $j$ that the label
of some representative of $\partial\Pi$ is~$r_j^{\pm1}$.
Clearly, the rank of a face is well-defined.
Let the rank of a face $\Pi$ be denoted by $\rank(\Pi)$.
Note that if $\rank(\Pi_1)\ge\rank(\Pi_2)$, then
$|\partial\Pi_1|\ge |\partial\Pi_2|$.
The ranks of two faces in a special S-diagram over
$\langle\,\mathfrak A\,\|\,\mathcal R\,\rangle$
are equal if and only if the lengths of the maximal selected subpaths
of their contours are equal (follows from \thetag{\ref{display:(2.2)}}).

\begin{proposition}
\label{proposition:(5.1)}
Every weakly reduced special S-diagram over\/
$\langle\,\mathfrak A\,\|\,\mathcal R\,\rangle$
satisfies the condition\/ $\mathcal B(\lambda_1,\lambda_2)$\textup.
\end{proposition}

\begin{proof}
Let $\Delta$ be a weakly reduced special S-diagram over
$\langle\,\mathfrak A\,\|\,\mathcal R\,\rangle$.
Clearly, $\Delta$ satisfies~$\mathcal B_0$.

If $\Pi$ is a face of $\Delta$, $s$ is the maximal selected subpath
of $\partial\Pi$, and $t$ is the path such that
$\langle st\rangle=\partial\Pi$,
then $\ell (s)=(x_1^{m_j}x_2^{m_j}\dots x_n^{m_j})^{\pm1}$ and
$\ell (t)=w_j^{\mp1}$ where $j=\rank{\Pi}$.
Since $|w_j|\le\lambda_1|r_j|$ (see \thetag{\ref{display:(2.4)}}),
have that $|s|\ge(1-\lambda_1)|\partial\Pi|$.
Hence, $\Delta$ satisfies $\mathcal B_1(\lambda_1)$.

If $u$ is a double-selected oriented arc,
then $\ell (u)$ is a subword of a word of the form $x_l^mx_{l+1}^m$ or
$x_{l+1}^{-m}x_l^{-m}$ (since $\Delta$ is weakly reduced).
Therefore, if such a $u$ is a subpath of the contour of a face $\Pi$,
then $|u|<(2/n)|\partial\Pi|=\lambda_2|\partial\Pi|$.
Hence, $\Delta$ satisfies $\mathcal B_2(\lambda_2)$.
\end{proof}

\begin{corollary.proposition}
\label{corollary.proposition:1.(5.1)}
If\/ $\Delta$ is a weakly reduced special S-diagram over\/
$\langle\,\mathfrak A\,\|\,\mathcal R\,\rangle$ with at most\/ $3$ contours\textup,
and\/ $S$ is the number of selected external edge of\/ $\Delta$\textup,
then
$$
S\ge(1-2\mu)\!\sum_{\Pi\in\Delta(2)}\!|\partial\Pi|.
$$
\end{corollary.proposition}

\begin{proof}
Follows from the proposition and the Main Lemma.
\end{proof}

\begin{corollary.proposition}
\label{corollary.proposition:2.(5.1)}
The group presentation\/ $\langle\,\mathfrak A\,\|\,\mathcal R\,\rangle$ is
strongly aspherical in the sense that no spherical diagram over\/
$\langle\,\mathfrak A\,\|\,\mathcal R\,\rangle$ is weakly reduced\textup.
\end{corollary.proposition}

\begin{proof}
Every weakly reduced special spherical S-diagram over
$\langle\,\mathfrak A\,\|\,\mathcal R\,\rangle$
has at least $(1-2\mu)\cdot2$ selected external edges
and therefore does not exist ($\mu<1/2$).
\end{proof}

\begin{proposition}
\label{proposition:(5.2)}
Let\/ $\Delta$ be a special nondegenerate S-diagram over\/
$\langle\,\mathfrak A\,\|\,\mathcal R\,\rangle$\textup.
Let\/ $\mathfrak B$ be an arbitrary nonempty subset of\/ $\mathfrak A$\textup,
and let\/ $k=|\mathfrak B|$\textup.
Let\/ $S$ be the number of selected external edges of\/ $\Delta$
whose labels are in\/~$\mathfrak B$\textup.
Then
$$
S<\frac{k}{n}\sum_{\Pi\in\Delta(2)}\!|\partial\Pi|.
$$
\end{proposition}

\begin{proof}
If $\Pi$ is a face of $\Delta$, and $s$ is the maximal selected subpath
of $\partial\Pi$, then $\ell (s)=(x_1^{m_j}x_2^{m_j}\dots x_n^{m_j})^{\pm1}$
where $j=\rank{\Pi}$.
Therefore, the number of selected external edges of $\Delta$ incident to $\Pi$
whose labels are in $\mathfrak B$ is at most $(k/n)|s|<(k/n)|\partial\Pi|$.
Since this is true for every face $\Pi$ of $\Delta$,
the desired inequality follows.
\end{proof}

\begin{definition}
Let $\langle\,\mathfrak B\,\|\,\mathcal S\,\rangle$ be a finite group presentation.
A function $f\!:\mathbb N\to\mathbb R$ is called an \emph{isoperimetric function\/}
of $\langle\,\mathfrak B\,\|\,\mathcal S\,\rangle$ if
%and only if
for every group word $w$ over $\mathfrak B$ equal to $1$ modulo
the relations $r=1$, $r\in\mathcal S$, there exists a disc diagram $\Delta$ over
$\langle\,\mathfrak B\,\|\,\mathcal S\,\rangle$ with at most $f(|w|)$ faces
such that $w$ is the label of some representative of $\partial\Delta$.
The minimal isoperimetric function of $\langle\,\mathfrak B\,\|\,\mathcal S\,\rangle$
is called the \emph{Dehn function\/} of $\langle\,\mathfrak B\,\|\,\mathcal S\,\rangle$.
\end{definition}

\begin{proposition}
\label{proposition:(5.3)}
If\/ $\langle\,\mathfrak B\,\|\,\mathcal S\,\rangle$ is a finite subpresentation of\/
$\langle\,\mathfrak A\,\|\,\mathcal R\,\rangle$
\textup($\mathfrak B$ is a subset of\/ $\mathfrak A$\textup, and\/ $\mathcal S$ is
a subset of\/ $\mathcal R$\textup{),}
then the function\/ $f\!:\mathbb N\to\mathbb R$ given by the formula
$$
f(k)=\frac{k}{1-2\mu}
$$
is an isoperimetric function of\/ $\langle\,\mathfrak B\,\|\,\mathcal S\,\rangle$
\textup(recall that\/ $\mu<1/2$\textup{).}
\end{proposition}

\begin{proof}
If $\Delta$ is a weakly reduced disc diagram over
$\langle\,\mathfrak B\,\|\,\mathcal S\,\rangle$, and
$\langle q\rangle=\partial\Delta$, then,
as follows from Corollary~\ref{corollary.proposition:1.(5.1)}
of Proposition~\ref{proposition:(5.1)},
$$
|q|\ge(1-2\mu)\!\sum_{\Pi\in\Delta(2)}\!|\partial\Pi|\ge(1-2\mu)\|\Delta(2)\|.
$$
Hence, $\|\Delta(2)\|\le f(|q|)$.
\end{proof}

\begin{corollary.proposition}
\label{corollary.proposition:1.(5.3)}
Every finite subpresentation of\/ $\langle\,\mathfrak A\,\|\,\mathcal R\,\rangle$
presents a hyperbolic group\textup.
\end{corollary.proposition}

\begin{proof}
Use the characterization of hyperbolic groups in terms of
isoperimetric functions of their finite presentations.
According to Theorem~2.5, Theorem~2.12, and Corollary of the latter
in \cite{ABCFLMSS:1991:nwhg}, a group is hyperbolic if and only if it has
a finite presentation with a linear isoperimetric function.
(Note that in \cite{ABCFLMSS:1991:nwhg} all isoperimetric functions in the sense of
the last definition are called ``Dehn functions.'')
\end{proof}

Let \emph{deg-lex\/} be the order on the set of all group words
over $\mathfrak A$ described as follows.
Two group words are compared first by length,
and second alphabetically
according to the following order on the group letters:
$x_1<x_1^{-1}<x_2<\dots<x_n<x_n^{-1}$.
For example, $x_nx_n<x_1x_1x_1$ and $x_1x_2x_3<x_1x_3x_2$.

\begin{proposition}
\label{proposition:(5.4)}
If the order on reduced group words used in
Section\/~\textup{\ref{section:group_construction}}
for choosing\/ $\{w_i\}_{i=1,\dots}$ is deg-lex\textup,
$N$ is a positive integer\textup,
and\/ $m_i=N|w_i|+i$ for every\/ $i$ for which\/ $r_i$ is defined\textup,
then the presentation\/ $\langle\,\mathfrak A\,\|\,\mathcal R\,\rangle$
is recursive\textup.
\end{proposition}

\begin{proof}
If the presentation $\langle\,\mathfrak A\,\|\,\mathcal R\,\rangle$
is finite, then it is recursive.
(It will be shown in Section~\ref{section:comments}
that it cannot be finite.)
Now, assume that the presentation is infinite.

Let $C$ be the set of all $4$-tuples $(\mathcal S,E,u,v)$ such that
$\mathcal S$ is a finite set of group words over $\mathfrak A$,
$E$ is a rational number, $u$ and $v$ are group words over $\mathfrak A$, and
there exists a disc diagram over
$\langle\,\mathfrak A\,\|\,\mathcal S\,\rangle$
with at most $E$ edges such that the label of some representative of its
contour is~$uv^{-1}$.
The set $C$ is recursive.

Let $q$ be a rational number such that $q\ge1/(1-2\mu)$.
Let $D$ be the set of all $3$-tuples $(\mathcal S,u,v)$ such that
$$
\Bigl(\mathcal S,\frac{1+qL}{2}(|u|+|v|),u,v\Bigr)\in C
$$
if $L$ is the maximum of the lengths of the elements of~$\mathcal S$.
The set $D$ is recursive.
This set describes consequences of a given finite set of relations
as explained below.

Let $f$ be the function $\mathbb N\to\mathbb R,k\mapsto qk$.
Consider an arbitrary finite $\mathcal S$ such that $f$ is
an isoperimetric function
of $\langle\,\mathfrak A\,\|\,\mathcal S\,\rangle$.
Let $u$ and $v$ be arbitrary group words over~$\mathfrak A$.
Then the relation $u=v$ is a consequence of the relations
of $\langle\,\mathfrak A\,\|\,\mathcal S\,\rangle$ if and only if $(\mathcal S,u,v)\in D$.

Consider the following algorithm:

\begin{quotation}
\raggedright
\tt
\begin{description}
\item[Input]
	a positive integer $j$.
\item[Step~1]
	Produce the set $\mathcal S$ of outputs of this algorithm
	on the inputs $1$, \dots, $j-1$.
	(Do not stop if the algorithm does not stop on at leas one of those inputs;
	if $j=1$, then $\mathcal S=\emptyset$.)
	Let $L$ be the maximum of the lengths of the elements of $\mathcal S$
	if $\mathcal S$ is nonempty, and $0$ otherwise.
\item[Step~2]
	Find the least (with respect to deg-lex) group word $w$ over $\mathfrak A$
	such that for every regular group word $u$ of length at most $(n+1)|w|+n^4L$,
	$(\mathcal S,u,w)\not\in D$.
	(Do not stop if there is no such $w$.)
	Let $m=N|w|+j$.
\item[Output]
	$x_1^mx_2^m\dots x_n^mw^{-1}$.
\end{description}
\end{quotation}

It shall be shown by induction that on every input $j\in\mathbb N$,
the algorithm stops and gives $r_j$ as the output.
It will follow that $\mathcal R$ is recursive
because the length of $r_j$ is an increasing function of~$j$.

Take an arbitrary $k\in\mathbb N$.
If $k>1$, assume that it is already proved that
for every positive integer $j<k$,
the algorithm gives $r_j$ as the output on the input~$j$.
Consider the work of the algorithm on the input~$k$.

By the inductive assumption, the set $\mathcal S$ produced in Step~1
coincides with~$\mathcal R_{k-1}$.
In particular, $f$ is an isoperimetric function
of $\langle\,\mathfrak A\,\|\,\mathcal S\,\rangle$.
Note that $L=|r_{k-1}|$ if $k>1$, and $L=0$ if $k=1$.

By the choice of $w_k$, this word is not equal to any regular group word
modulo the relations of $\langle\,\mathfrak A\,\|\,\mathcal S\,\rangle$.
Therefore, Step~2 of the algorithm is carried out in finite time,
and the obtained group word $w$ is not greater than $w_k$ with respect to
the deg-lex order.

To complete the inductive step and the proof of this proposition,
it suffices to show that the group word $w$ obtained in Step~2 is~$w_k$.
Suppose $w$ is distinct from~$w_k$.

Since $w<w_k$ relative to deg-lex, it follows from the choice of $w_k$ that
$w$ equals some regular group word modulo the relations
of $\langle\,\mathfrak A\,\|\,\mathcal S\,\rangle$.
Let $u$ be such a regular group word.
By the choice of $w$,
$$
|u|>(n+1)|w|+n^4L.
$$

Let $\Delta$ be a weakly reduced special disc S-diagram
over $\langle\,\mathfrak A\,\|\,\mathcal S\,\rangle$
such that the label of some representative of $\partial\Delta$ is $uw^{-1}$.
Let $b$ and $p$ be the paths such that $\ell (b)=u$, $\ell (p)=w^{-1}$, and
$\langle bp\rangle=\partial\Delta$.
Let $b_1$, $b_2$, \dots, $b_n$ be the subpaths of $b$ such that
$b=b_1b_2\dots b_n$, and for every $i=1,2,\dots,n$, the label of $b_i$
is a power of~$x_i$.
For every $i=1,2,\dots,n$, let $B_i$ be the set of all edges that
lie on~$b_i$.
Let $\Sigma$ be the sum of the degrees of all the faces of~$\Delta$.
Let $T$ be the number of selected external edges of $\Delta$ that lie on~$b$.

Since $\ell (b)$ is regular, no two oriented edges of $b$ are inverse
to each other.
Hence,
$$
|b|\le\Sigma+|p|.
$$

Let $\Sigma_1$ be the sum of the degrees of all faces of $\Delta$ that
are incident with edges from at least $3$ distinct sets
from among $B_1$, \dots,~$B_n$.
If $1\le i_1<i_2<i_3\le n$, then there is at most one face that is incident to
edges from all three sets $B_{i_1}$, $B_{i_2}$,~$B_{i_3}$.
Therefore, $\Sigma_1\le n^3L$
(the degree of every face of $\Delta$ is at most~$L$).
Let $\Sigma_2$ be the sum of the degrees of all faces of $\Delta$
that are incident with edges from no more than $2$ distinct sets
from among $B_1$, \dots,~$B_n$.

One one hand,
$$
T\le\Sigma_1+\frac{2}{n}\Sigma_2\le n^3L+\frac{2}{n}\Sigma.
$$
On the other hand, by Proposition~\ref{proposition:(5.1)} and
the Main Lemma,
$$
T+|p|\ge(1-2\mu)\Sigma.
$$
Therefore,
$$
\aligned
(1-2\mu)\Sigma&\le n^3L+\frac{2}{n}\Sigma+|p|,\\
\Bigl(1-2\mu-\frac{2}{n}\Bigr)\Sigma &\le |p|+n^3L,\\
\frac{n-23}{n}\Sigma &\le |p|+n^3L,\\
\Sigma &\le n|p|+n^4L,\\
|b|&\le (n+1)|p|+n^4L.
\endaligned
$$
In other terms, $|u|\le (n+1)|w|+n^4L$.
This gives a contradiction.

Thus, $w=w_k$, and the inductive step is done.
\end{proof}

Everything is now ready to start proving properties of~$G$.

\begin{property}
\label{property:1}
For every\/ $(n-21)$-element subset\/ $I$ of the set\/ $\{1,\dots,n\}$\textup,
the system\/ $\{a_i\}_{i\in I}$ freely generates 
a free subgroup\/ $F$ of\/ $G$
\textup(of rank\/ $n-21$\textup{).}
Moreover\textup, every two elements in such a free subgroup\/ $F$
are conjugate in\/ $G$ only if they are conjugate in\/~$F$\textup.
%\textup(\ie\textup, $F$ is Frattini embedded into\/ $G$ by
%the inclusion monomorphism\/\textup{).}
\end{property}

\begin{proof}
Suppose that Property~\ref{property:1} does not hold.
Take a set $I\subset\{1,\dots,n\}$ of $n-21$ elements
that provides a counterexample.
Let $\mathfrak B=\{\,x_i\mid i\in I\,\}$.

If the system $\{a_i\}_{i\in I}$ does not freely generate
a free subgroup of $G$,
then there is a weakly reduced disc diagram $\Delta$ over
$\langle\,\mathfrak A\,\|\,\mathcal R\,\rangle$ such that
the label of some (therefore every) representative of its contour
is a cyclically reduced nonempty group word
over the alphabet~$\mathfrak B$.
Clearly, such a diagram $\Delta$ is nondegenerate (has a face).

If the system $\{a_i\}_{i\in I}$ does freely generate
a free subgroup $F$ of $G$,
but there are two elements of $F$ that are conjugate in $G$ but not in $F$,
then let $v_1$ and $v_2$ be two cyclically reduced group words
over $\mathfrak B$ that represent two such elements.
The group words $v_1$ and $v_2$ are not cyclic shifts of each other,
and neither of them represents the identity of~$G$.
Since $[v_1]$ and $[v_2]$ are nontrivial but are conjugate in $G$,
there exists a weakly reduced annular diagram $\Delta$ over
$\langle\,\mathfrak A\,\|\,\mathcal R\,\rangle$ such that
the label of some representative of one of its contours is $v_1$,
and the label of some representative of its other contour is $v_2^{-1}$
(see Lemma V.5.2 in \cite{LyndonSchupp:2001:cgt}).
Clearly, such $\Delta$ is nondegenerate.

Thus, to obtain a contradiction and complete the proof, it is enough
to show that there is no weakly reduced nondegenerate
disc nor annular diagram over $\langle\,\mathfrak A\,\|\,\mathcal R\,\rangle$
with fewer than $n-20$ distinct basic letters on its contour(s).

Let $\Delta$ be an arbitrary special weakly reduced nondegenerate
disc or annular S-diagram over
$\langle\,\mathfrak A\,\|\,\mathcal R\,\rangle$.
(Recall that a special selection exists on every diagram over
$\langle\,\mathfrak A\,\|\,\mathcal R\,\rangle$.)
Suppose that $\Delta$ has fewer than $n-20$ distinct basic letters
of $\mathfrak A$ on its contour(s).
Let $S$ be the total number of selected external edges of~$\Delta$.
On one hand, by Proposition~\ref{proposition:(5.2)},
$$
S<\frac{n-21}{n}\!\sum_{\Pi\in\Delta(2)}\!|\partial\Pi|.
$$
On the other hand, by Corollary~\ref{corollary.proposition:1.(5.1)}
of Proposition~\ref{proposition:(5.1)},
$$
S\ge(1-2\mu)\!\sum_{\Pi\in\Delta(2)}\!|\partial\Pi|.
$$
Since
$$
1-2\mu>\frac{n-21}{n},
$$
one has a contradiction.
\end{proof}

The proof of the next property is not so straightforward.

\begin{property}
\label{property:2}
For every element\/ $g$ of\/ $G$\textup, there exist unique\/
$k_1,\dots,k_n\in\mathbb Z$ such that\/ $g=a_1^{k_1}\dots a_n^{k_n}$\textup.
\end{property}

%% AMS-LaTeX
%% Paper 2 - Version GAMMA
%% Subpart of Part 05
%% (The hard part of the paper)
%%%%%%%%%%%%%%%%%%%%%%%%%%%%%%%%%%%%%%%%%%%%%%%%%%%%%%%%%%%%%%%%%%%%%%%%%%%

%!TEX root = paper2-2005.tex

\par
The uniqueness is the ``hard part'' of Property~\ref{property:2}.
It shall be proved by contradiction.
Main steps of the proof are stated below as Lemmas
\ref{lemma:(5.1)}--\ref{lemma:(5.16)}.
These lemmas share some common assumptions about an S-diagram~$\Delta$.

Assume that $\Delta$ is a special disc S-diagram over
$\langle\,\mathfrak A\,\|\,\mathcal R\,\rangle$
whose contour has a representative of the form $p_1p_2$ such that
$\ell (p_1)$ and $\ell (p_2^{-1})$ are distinct regular group words,
and $\ell (p_1p_2)$ is cyclically reduced.
Suppose, moreover, that $\Delta$ is such an S-diagram with the minimal
possible number of faces.
In particular, if $\Delta'$ is a nontrivial disc diagram over
$\langle\,\mathfrak A\,\|\,\mathcal R\,\rangle$,
$p$ is a representative of $\partial\Delta'$, and
$\ell (p)$ is regular or counter-regular,
then $\Delta'$ has no fewer faces than~$\Delta$.

\begin{lemma}
\label{lemma:(5.1)}
The diagram\/ $\Delta$ is weakly reduced\textup.
\end{lemma}

\begin{proof}
This easily follows from the minimality of the number of faces of~$\Delta$.
\end{proof}

\begin{lemma}
\label{lemma:(5.2)}
The diagram\/ $\Delta$ is simple\textup.
\end{lemma}

\begin{proof}
Use the minimality of $\Delta$ again.
Some maximal simple disc S-subdiagram of $\Delta$ satisfies all
the assumptions made about $\Delta$
(this follows, for example, from Proposition~3.1
of~\cite{Muranov:2005:dsmcbgbsg}),
and hence it must be the whole of~$\Delta$.
\end{proof}

\begin{lemma}
\label{lemma:(5.3)}
If\/ $\Delta'$ is a disc S-subdiagram of\/ $\Delta$\textup,
and\/ $S$ is the number of selected external edge of\/ $\Delta'$
\textup(external in\/ $\Delta'$\textup{),}
then
$$
S\ge(1-2\mu)\!\sum_{\Pi\in\Delta'(2)}\!|\partial\Pi|.
$$
\end{lemma}

\begin{proof}
Follows directly from Lemma~\ref{lemma:(5.1)} and
Corollary~\ref{corollary.proposition:1.(5.1)} of
Proposition~\ref{proposition:(5.1)}.
\end{proof}

\begin{lemma}
\label{lemma:(5.4)}
If\/ $\Delta'$ is a disc subdiagram of\/ $\Delta$\textup,
$p_1$ and\/ $p_2$ are paths\textup,
$\langle p_1p_2\rangle=\partial\Delta'$\textup,
and the reduced forms of\/ $\ell (p_1)$
and\/ $\ell (p_2^{-1})$ are regular\textup,
then either\/ $\Delta'$ is degenerate\textup,
or it coincides with\/~$\Delta$\textup.
\end{lemma}

\begin{proof}
Suppose $\Delta'$, $p_1$, and $p_2$ are such as in the hypotheses of
the lemma.

Case~1: $\ell (p_1)$ and $\ell (p_2^{-1})$ are freely equal.
Then the label of every representative of $\partial\Delta'$
is freely trivial
and therefore $\Delta'$ is degenerate.
Indeed, if $\Delta'$ was nondegenerate,
there would exist a disc diagram over
$\langle\,\mathfrak A\,\|\,\mathcal R\,\rangle$
with the same label of the contour as $\Delta$ but with fewer faces
(faces of $\Delta'$, possibly together with some of the others,
could be ``eliminated'' from $\Delta$),
which would contradict the minimality of~$\Delta$.

Case~2: $\ell (p_1)$ and $\ell (p_2^{-1})$ are not freely equal.
Then take $\Delta'$ and repeatedly \emph{fold\/} labelled external edges
and ``cut off,'' whenever necessary, ``branches'' with freely trivial
contour labels until one obtains a disc diagram $\Delta''$ such that
the labels of the representatives of $\partial\Delta''$
are cyclically reduced.
Clearly, $\Delta''$ is a nondegenerate diagram over
$\langle\,\mathfrak A\,\|\,\mathcal R\,\rangle$,
the label of every representative of $\partial\Delta''$ is
a cyclically reduced form of $\ell(p_1p_2)$,
and $\|\Delta''(2)\|\le\|\Delta'(2)\|$.
By the minimality of $\Delta$, such a diagram $\Delta''$ cannot have
fewer faces than~$\Delta$.
Therefore, $\Delta'(2)=\Delta(2)$ and,
since $\partial\Delta$ is reduced,
$\Delta'$ coincides with~$\Delta$.
\end{proof}

\begin{lemma}
\label{lemma:(5.5)}
The diagram\/ $\Delta$ has more than\/ $1$ face\textup.
\end{lemma}

\begin{proof}
This easily follows from the form of elements of~$\mathcal R$.
\end{proof}

\begin{lemma}
\label{lemma:(5.6)}
There are at least\/ $n-20$ distinct basic letters
on\/~$\partial\Delta$\textup.
\end{lemma}

\begin{proof}
Follows from Property~\ref{property:1}.
\end{proof}

\begin{figure}
\includegraphics{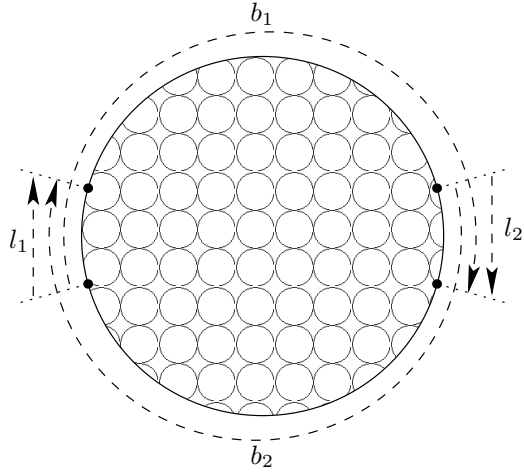}
\caption{The paths $b_1$, $b_2$, $l_1$, and $l_2$ in $\Delta$.}
\label{figure:01}
\end{figure}

Let $b_1$ be the maximal subpath of $\partial\Delta$ whose label is
regular but not counter-regular (is not a letter power).
Let $b_2$ be the maximal subpath of $\partial\Delta$ whose label is
counter-regular but not regular.
Since, according to Lemma~\ref{lemma:(5.6)},
there are at least $2$ distinct basic
letters on $\partial\Delta$, the paths $b_1$ and $b_2$ are well-defined.
(If, say, the label of some representative of $\partial\Delta$ is
$x_1x_2x_3\dots x_{42}x_{43}x_{42}x_2x_1$, then
$\ell (b_1)=x_1^2x_2\dots x_{43}$, and $\ell (b_2)=x_{43}x_{42}x_2x_1^2$.)
The label of each of the paths $b_1$ and $b_2$ has
at least $2$ distinct basic letters.
Let $l_1$ be the initial subpath of $b_1$ that is
a terminal subpath of $b_2$,
and let $l_2$ be the initial subpath of $b_2$ that is
a terminal subpath of~$b_1$.
Both $l_1$ and $l_2$ are nontrivial (see Fig.~\ref{figure:01}).
The labels of $l_1$ and $l_2$ are letter powers.

\begin{lemma}
\label{lemma:(5.7)}
Every selected external arc of\/ $\Delta$ lies
on at least one of the paths\/ $b_1$ or\/~$b_2$\textup.
\end{lemma}

\begin{proof}
Follows from the fact that the label of every selected oriented arc
is either regular or counter-regular.
%Consider an arbitrary selected external arc $u$ of~$\Delta$.
%Then each of the two corresponding oriented arcs has either
%regular or counter-regular label.
%Consider the oriented arc corresponding to $u$ that is a subpath of
%$\partial\Delta$.
%Clearly, if $\ell (u)$ is simultaneously regular and counter-regular
%(\ie, is a letter power),
%then $u$ lies on $b_1$ or $b_2$ (possibly on both);
%if $\ell (u)$ is regular but not counter-regular, then $u$ lies on $b_1$;
%if $\ell (u)$ is counter-regular but not regular, then $u$ lies on~$b_2$.
\end{proof}

The following observation is obvious but deserves mentioning
because it is used implicitly in the proofs of Lemmas
\ref{lemma:(5.8)}, \ref{lemma:(5.10)},
and \ref{lemma:(5.13)} several times.

Suppose $p_1$, $p_2$, and $q$ are paths in a diagram over
$\langle\,\mathfrak A\,\|\,\mathcal R\,\rangle$, and $q$ is nontrivial.
Then
\begin{itemize}
\item
	if the products $p_1q$ and $q^{-1}p_2$ are defined
	and their labels are regular,
	then $\ell (q)$ is a letter power, and
	the reduced form of $\ell (p_1p_2)$ is regular;
\item
	if the products $p_1q$ and $qp_2$ are defined
	and their labels are regular,
	then $\ell (p_1qp_2)$ is regular.
\end{itemize}

Since Lemma \ref{lemma:(5.7)} and many of the remaining lemmas deal with
selected external arcs or selected external oriented arcs,
it is advisable to review these concepts.
An oriented arc of a graph (or a $2$-complex) is a simple path all of whose
vertices, except, possibly, the end-vertices, have degree two
in the given graph.
A (non-oriented) arc is a pair of mutually inverse oriented arcs.
An arc is incident to a face $\Pi$ if (and only if) at least one of
the associated oriented arcs is a subpath of~$\partial\Pi$.
Every external arc of a (semi-)simple map is incident to exactly one face.
An oriented arc is selected if (and only if) it is selected as a path
(in particular, it must be a subpath of the contour of some face).
An external arc is selected if (and only if) it is incident to some
face $\Pi$, and the associated oriented arc that is a subpath
of $\partial\Pi$ is selected.
Distinct maximal selected external arcs of $\Delta$ never overlap.

\begin{lemma}
\label{lemma:(5.8)}
In the diagram\/ $\Delta$\textup, every face is incident with at most\/
$1$ maximal selected external arc\textup.
\end{lemma}

\begin{proof}
Consider an arbitrary face $\Pi$ of~$\Delta$.
Suppose $\Pi$ is incident with at least $2$ distinct maximal selected
external arcs.

Let $u_1$ and $u_2$ be distinct maximal selected external
oriented arcs of $\Delta$ that are subpaths of $\partial\Pi$.
(They do not overlap.)
Let $s$ be the maximal selected subpath of $\partial\Pi$.
Without loss of generality, assume that $u_1$ precedes $u_2$
as a subpath of~$s$.
Let $v$ be such a path that $u_1vu_2$ is a subpath of~$s$.
Since $\Delta$ is simple (see Lemma~\ref{lemma:(5.2)}),
and $u_1$ and $u_2$ are maximal, it follows that $v$ is nontrivial, and
the first and last oriented edges of $v$ are internal in~$\Delta$.

Suppose $v$ has a simple cyclic subpath~$p$.
Then $\langle p\rangle$ is the contour
of a proper simple disc subdiagram of~$\Delta$.
Since $\ell (p)$ is regular or counter-regular,
this contradicts the minimality of~$\Delta$.
Hence, the path $v$ is simple and not cyclic
(equivalently, every vertex occurs in $v$ at most once).

Let $v_1$ be the minimal subpath of $\partial\Delta$
such that $u_1^{-1}$ and $u_2^{-1}$ are respectively
its initial and terminal subpaths.
Let $v_2$ be the path
such that $\langle v_1v_2\rangle=\partial\Delta$.
The paths $v_1v^{-1}$ and $v_2v$ are cyclically reduced and represent
the contours of two disc subdiagrams of~$\Delta$.
Let $\Delta_1'$ and $\Delta_2'$ be the subdiagrams
with the contours $\langle v_1v^{-1}\rangle$ and $\langle v_2v\rangle$,
respectively.
Note that $\Pi\in\Delta_1'(2)$, and $\Delta(2)$ is the disjoint union of
$\Delta_1'(2)$ and $\Delta_2'(2)$.
The diagram $\Delta_1'$ is simple.
The diagram $\Delta_2'$ is nondegenerate.
Both $\Delta_1'$ and $\Delta_2'$ are proper subdiagrams of~$\Delta$.

The labels of $u_1$ and $u_2$ are both regular or both counter-regular.
Therefore, each of the paths $u_1^{-1}$ and $u_2^{-1}$ is a subpath
of one of the paths $b_1$ or $b_2$ (see Lemma~\ref{lemma:(5.7)}).
Consider the following 4 cases
(see Fig.~\ref{figure:02};
shaded subdiagrams are the ``sources of contradiction''):

\begin{figure}
\includegraphics{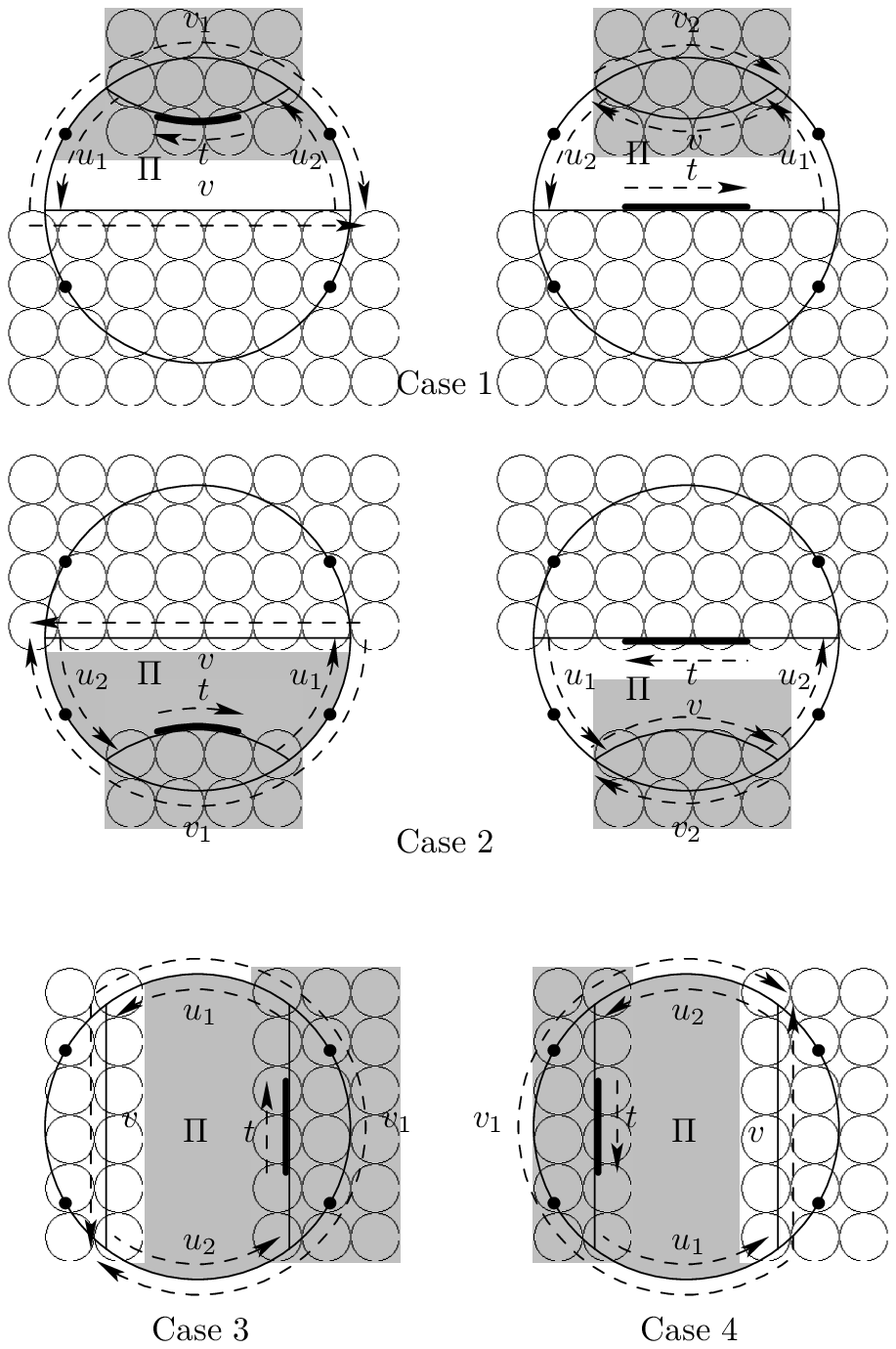}
\caption{Cases 1--4, Lemma \ref{lemma:(5.8)}.}
\label{figure:02}
\end{figure}

Case~1: $u_1^{-1}$ and $u_2^{-1}$ are subpaths of~$b_1$.
Then either $v_1$ or $u_2^{-1}v_2u_1^{-1}$ is a subpath of $b_1$ as well.
Suppose $v_1$ is a subpath of~$b_1$.
Then $\ell (v_1)$ is regular,
$\ell (u_1)$ and $\ell (u_2)$ are letter powers, and
$\ell (u_1vu_2)$ is regular.
Therefore, $\langle v_1v^{-1}\rangle$ cannot be the contour
of a proper nondegenerate disc subdiagram (see Lemma~\ref{lemma:(5.4)}).
This gives a contradiction.
Suppose $u_2^{-1}v_2u_1^{-1}$ is a subpath of~$b_1$.
Then $\ell (u_2^{-1}v_2u_1^{-1})$ is regular, and $\ell (u_1vu_2)$
is counter-regular.
Therefore, $\langle v_2v\rangle$ cannot be the contour
of a proper nondegenerate disc subdiagram.
This gives a contradiction.

Case~2: $u_1^{-1}$ and $u_2^{-1}$ are subpaths of~$b_2$.
%Then either $v_1$ or $u_2^{-1}v_2u_1^{-1}$ is a subpath of $b_2$ as well.
%Suppose $v_1$ is a subpath of~$b_2$.
%Then $\ell (v_1)$ is counter-regular,
%$\ell (u_1)$ and $\ell (u_2)$ are letter powers,
%and $\ell (u_1vu_2)$ is counter-regular.
%Therefore, $\langle v_1v^{-1}\rangle$ cannot be the contour
%of a proper nondegenerate disc subdiagram.
%This gives a contradiction.
%Suppose $u_2^{-1}v_2u_1^{-1}$ is a subpath of~$b_1$.
%Then $\ell (u_2^{-1}v_2u_1^{-1})$ is counter-regular, and $\ell (u_1vu_2)$
%is regular.
%Therefore, $\langle v_2v\rangle$ cannot be the contour
%of a proper nondegenerate disc subdiagram.
%This gives a contradiction.
%Note that a contradiction in this case may also be derived from
%a contradiction in Case~1 by passing to the mirror copy of~$\Delta$.
When proving impossibility of Case~1, it was essentially shown that
no face of $\Delta$ can be incident with $2$ distinct maximal selected
external arcs lying on~$b_1$.
Since all the assumptions made about $\Delta$ hold for its
mirror copy as well, the same statement appropriately reformulated
must hold for the mirror copy of~$\Delta$.
Namely, no face of the mirror copy of $\Delta$ can be incident with
$2$ distinct maximal selected external arcs lying on $b_2^{-1}$
($b_2^{-1}$ plays the same role for the mirror copy of $\Delta$ as
$b_1$ does for $\Delta$).
This means that Case~2 is impossible.

Case~3: $u_1^{-1}$ is a subpath of $b_1$, and
$u_2^{-1}$ is a subpath of~$b_2$.
Pick a subpath $p_1$ of $b_1$ and a subpath $p_2$ of $b_2$ such that
$u_1^{-1}$ is an initial subpath of $p_1$,
$u_2^{-1}$ is a terminal subpath of $p_2$, and
$p_1p_2=v_1$.
Then $\ell (p_1)$ is regular,
$\ell (p_2)$ is counter-regular, and
$\ell (u_2^{-1}v^{-1}u_1^{-1})$ is regular or counter-regular.
If $\ell (u_2^{-1}v^{-1}u_1^{-1})$ is regular,
then $\ell (v^{-1}p_1)$ is regular.
If $\ell (u_2^{-1}v^{-1}u_1^{-1})$ is counter-regular,
then $\ell (p_2v^{-1})$ is counter-regular.
In either case $\langle p_2v^{-1}p_1\rangle=\langle v_1v^{-1}\rangle$
cannot be the contour of a proper nondegenerate disc subdiagram.
This gives a contradiction.

Case~4: $u_1^{-1}$ is a subpath of $b_2$, and
$u_2^{-1}$ is a subpath of~$b_1$.
Pick a subpath $p_1$ of $b_2$ and a subpath $p_2$ of $b_1$ such that
$u_1^{-1}$ is an initial subpath of $p_1$,
$u_2^{-1}$ is a terminal subpath of $p_2$, and
$p_1p_2=v_1$.
Then $\ell (p_1)$ is counter-regular,
$\ell (p_2)$ is regular, and
$\ell (u_2^{-1}v^{-1}u_1^{-1})$ is regular or counter-regular.
If $\ell (u_2^{-1}v^{-1}u_1^{-1})$ is regular,
then $\ell (p_2v^{-1})$ is regular.
If $\ell (u_2^{-1}v^{-1}u_1^{-1})$ is counter-regular,
then $\ell (v^{-1}p_1)$ is counter-regular.
In either case $\langle p_2v^{-1}p_1\rangle=\langle v_1v^{-1}\rangle$
cannot be the contour of a proper nondegenerate disc subdiagram.
This gives a contradiction.

(Case~4 is symmetric with Case~3, but this symmetry is not simply
the mirror symmetry as between Cases 1 and~2.
Say, if Case~4 was indeed the case for $\Delta$, then
the mirror copy of $\Delta$ would be just another example for
the same Case~4, not Case~3.)

A contradiction is obtained in each of the 4 cases,
and no other cases exist.
%Thus, the face $\Pi$ is incident with at most $1$ maximal selected
%external arc of~$\Delta$.
\end{proof}

\begin{lemma}
\label{lemma:(5.9)}
There are a face\/ $\Pi$ and a path\/ $s'$ in\/ $\Delta$
such that\/\textup:
\begin{itemize}
\item
	$s'$ is a selected subpath of\/ $\partial\Pi$\textup;
\item
	$s'$ is a maximal selected external oriented arcs of\/
	$\Delta$\textup;
\item
	the\/ $n-40$ basic letters\/
	$x_{21}$\textup, $x_{22}$\textup, \dots\textup, $x_{n-20}$
	all occur in\/~$\ell (s')$\textup.
\end{itemize}
\end{lemma}

\begin{proof}
For every face $\Pi$ of $\Delta$,
let $S(\Pi)$ denote the number of selected external edges of $\Delta$
incident to~$\Pi$.
By Lemma~\ref{lemma:(5.3)},
$$
\sum_{\Pi\in\Delta(2)}\!S(\Pi)
\ge(1-2\mu)\!\sum_{\Pi\in\Delta(2)}\!|\partial\Pi|.
$$
Hence, there exists a face $\Pi$ of $\Delta$ such that
$$
S(\Pi)\ge(1-2\mu)|\partial\Pi|>\frac{n-21}{n}|\partial\Pi|
$$
(recall that $1-2\mu>1-21/n$).

Let $\Pi$ be such a face as above.
Let $s$ be the maximal selected subpath of $\partial\Pi$.
Let $s'$ be the (only) maximal selected external oriented arc of $\Delta$
that is a subpath of $s$
(see Lemma~\ref{lemma:(5.8)}).
Then
$$
|s'|=S(\Pi)>\frac{n-21}{n}|\partial\Pi|>\frac{n-21}{n}|s|.
$$
Let $j=\rank(\Pi)$.
Then either $\ell (s)=x_1^{m_j}x_2^{m_j}\dots x_n^{m_j}$
or $\ell (s)=x_n^{-m_j}x_{n-1}^{-m_j}\dots x_1^{-m_j}$.
Therefore, $\ell (s')$ has at least $n-20$ distinct basic letters, and
all of the basic letters $x_{21}$, $x_{22}$, \dots, $x_{n-20}$
occur in it.
\end{proof}

Lemmas \ref{lemma:(5.1)}--\ref{lemma:(5.9)},
as well as \ref{lemma:(5.10)}--\ref{lemma:(5.16)},
assert some properties of~$\Delta$.
Observe that none of these properties in fact can distinguish
between $\Delta$ and its mirror copy, \ie,
each of these properties holds for $\Delta$ if and only if it holds
for the mirror copy of~$\Delta$.
Thus, at this point analogues of
Lemmas \ref{lemma:(5.1)}--\ref{lemma:(5.9)} for the mirror copy
of $\Delta$ shall be assumed proved.
Moreover, the initial assumptions about the S-diagram $\Delta$ are also
true about its mirror copy.

Even though the statement of Lemma~\ref{lemma:(5.7)} is about
the S-diagram $\Delta$ and the paths $b_1$ and $b_2$,
it still may be viewed as an assertion of a property of $\Delta$
because, according to the way they are chosen,
$b_1$ and $b_2$ are uniquely determined for the given~$\Delta$.
Therefore, the analog of Lemma~\ref{lemma:(5.7)}
for the mirror copy of $\Delta$ states:
\begin{quote}
Every selected external arc of the mirror copy of $\Delta$ lies
on at least one of the paths $b_2^{-1}$ or~$b_1^{-1}$.
\end{quote}

In the proofs of Lemmas \ref{lemma:(5.10)}--\ref{lemma:(5.13)},
it is convenient in some cases to pass to the mirror copy of $\Delta$
to reduce the number of cases to consider.

\begin{lemma}
\label{lemma:(5.10)}
Let\/ $\Pi_1$ be a face of\/ $\Delta$ incident to a selected
external edge of\/ $\Delta$\textup,
and let\/ $\Pi_2$ be another face of\/~$\Delta$\textup.
Let\/ $s_1$ and\/ $s_2$ be the maximal selected subpaths of\/
$\partial\Pi_1$ and\/ $\partial\Pi_2$\textup, respectively\textup.
Let\/ $s_1'$ be the maximal selected external oriented arc of\/
$\Delta$ that is a subpath of\/~$s_1$\textup.
Let\/ $s_{1-}'$ and\/ $s_{1+}'$ be the paths such that\/
$s_1=s_{1-}'s_1's_{1+}'$\textup.
Let\/ $q_1$ be the path such that\/
$\langle s_1^{\prime-1}q_1\rangle=\partial\Delta$\textup.
Suppose there are at least\/ $2$ distinct basic letters
in\/~$\ell (s_1')$\textup.
Suppose the paths\/ $s_2$ and\/ $(s_{1-}'q_1s_{1+}')^{-1}$
have a common oriented edge\textup.
Then they have exactly one maximal common nontrivial subpath\textup.
\end{lemma}

\begin{proof}
Observe that to prove that there is exactly one maximal common
nontrivial subpath of $s_2$ and $(s_{1-}'q_1s_{1+}')^{-1}$,
it suffices to prove the same for $s_2^{-1}$ and $s_{1-}'q_1s_{1+}'$.
This justifies passing to the mirror copy of~$\Delta$.

Since the paths $s_2$ and $(s_{1-}'q_1s_{1+}')^{-1}$ have a common
oriented edge, they have at least one maximal common nontrivial subpath.

\begin{figure}
\includegraphics{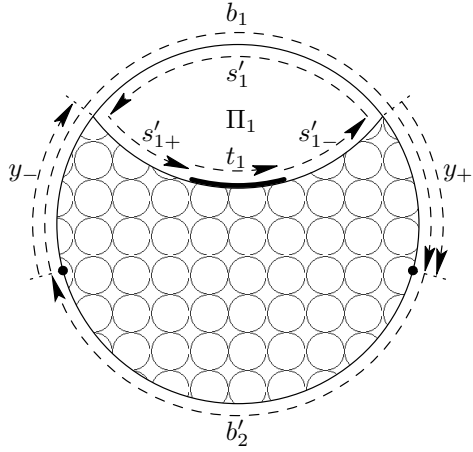}
\caption{The face $\Pi_1$ in $\Delta$, Lemma \ref{lemma:(5.10)}.}
\label{figure:03}
\end{figure}

The label of $s_1'$ is regular or counter-regular.
Therefore, $s_1'$ is a subpath of $b_1^{-1}$ or~$b_2^{-1}$.
If it is a subpath of $b_2^{-1}$ but not of $b_1^{-1}$, pass from
$\Delta$, $b_1$, $b_2$, $l_1$, $l_2$,
$s_1$, $s_2$, $s_{1-}'$, $s_{1+}'$, $s_1'$, $q_1$
to the mirror copy of $\Delta$,
$b_2^{-1}$, $b_1^{-1}$, $l_1^{-1}$, $l_2^{-1}$, $s_1^{-1}$, $s_2^{-1}$,
$s_{1+}^{\prime-1}$, $s_{1-}^{\prime-1}$, $s_1^{\prime-1}$, $q_1^{-1}$,
respectively.
Hence, it may and shall be assumed that $s_1'$ is a subpath of~$b_1^{-1}$
(see Fig.~\ref{figure:03}).

Let $t_1$ be the path such that $\langle s_1t_1\rangle=\partial\Pi_1$.
Let $y_-$ and $y_+$ be the paths such that $b_1=y_-s_1^{\prime-1}y_+$.
Let $b_2'$ be the subpath of $b_2$ such that
$\langle b_1b_2'\rangle=\partial\Delta$.
Then $b_2=l_2b_2'l_1$ and $q_1=y_+b_2'y_-$.

Since $\ell (b_1^{-1})$ is counter-regular, and $\ell (s_1')$ is not
a letter power (it has at least $2$ distinct basic letters),
$\ell (s_1)$ is counter-regular as well.
Hence, $\ell (s_{1+}^{\prime-1}s_1^{\prime-1}y_+)$ and
$\ell (y_-s_1^{\prime-1}s_{1-}^{\prime-1})$ are regular.

Let $j_1=\rank(\Pi_1)$.
Then $\ell (s_1)=x_n^{-m_{j_1}}x_{n-1}^{-m_{j_1}}\dots x_1^{-m_{j_1}}$
and $\ell (t_1)=w_{j_1}$.

Since $\ell (s_1')$ has at least $2$ distinct basic letters,
$\ell (y_-s_{1+}')$ and $\ell (s_{1-}'y_+)$ have disjoint sets of
basic letters.
This also implies that $l_1$ is a proper initial subpath of
$y_-s_1^{\prime-1}$,
and $l_2$ is a proper terminal subpath of $s_1^{\prime-1}y_+$.

Suppose $s_2$ and $(s_{1-}'q_1s_{1+}')^{-1}$ have at least two distinct
maximal common nontrivial subpaths.
Then let $u_1$ and $u_2$ be two paths such that:
\begin{itemize}
\item
	$u_1$ and $u_2$ are nontrivial subpaths of distinct maximal common
	subpaths of $s_2$ and $(s_{1-}'q_1s_{1+}')^{-1}$,
\item
	$u_1$ precedes $u_2$ as a subpath of $s_2$, and
\item
	each of the paths $u_1$ and $u_2$ is a maximal common subpath
	of $s_2$ and one of the paths $y_-^{-1}$, $y_+^{-1}$,
	$s_{1-}^{\prime-1}$, $s_{1+}^{\prime-1}$, or~$b_2^{\prime-1}$.
\end{itemize}

Clearly, such $u_1$ and $u_2$ exist and are non-overlapping oriented arcs
of~$\Delta$.
Let $v$ be the path such that $u_1vu_2$ is a subpath of~$s_2$.
Since $u_1$ and $u_2$ are subpaths of distinct maximal common subpaths
of $s_2$ and $(s_{1-}'q_1s_{1+}')^{-1}$,
the path $u_1vu_2$ is not a subpath of $(s_{1-}'q_1s_{1+}')^{-1}$.

Consider the following 16 cases
(see Fig.~\ref{figure:04}--\ref{figure:08}):

(Notice that in a certain sense Cases 6, 8, 10, 12, 14, and 16 are
``symmetric'' with Cases 5, 7, 9, 11, 13, and 15, respectively.)

Case~1: each of the paths $u_1$ and $u_2$ is a subpath of~$q_1^{-1}$.
Then $u_1$ and $u_2$ are external selected oriented arcs of~$\Delta$.
Since there is only one maximal selected external arc of $\Delta$
incident to $\Pi_2$ (see Lemma~\ref{lemma:(5.8)}),
the path $u_1vu_2$ is a subpath of~$q_1^{-1}$.
This gives a contradiction.

\begin{figure}
\includegraphics{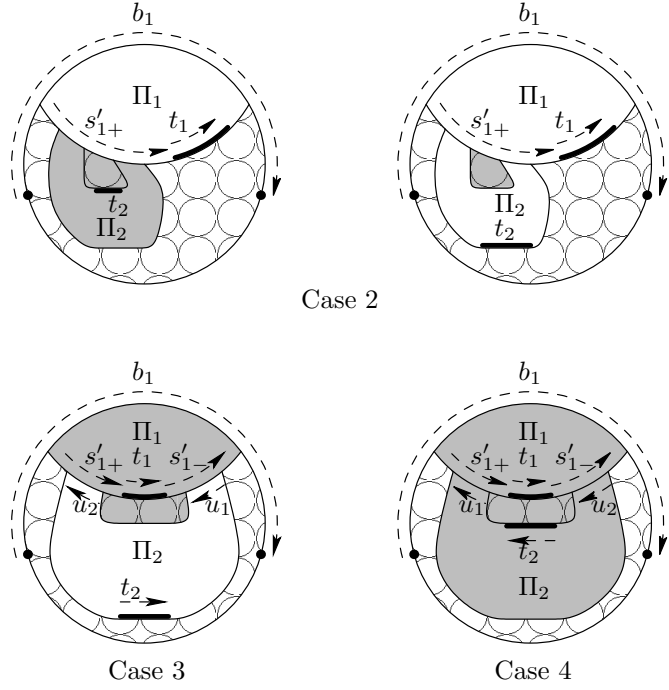}
\caption{Cases 2, 3, 4, Lemma \ref{lemma:(5.10)}.}
\label{figure:04}
\end{figure}

Case~2: both $u_1$ and $u_2$ are subpaths of one of the paths
$s_{1-}^{\prime-1}$ or~$s_{1+}^{\prime-1}$.
Let $z$ be one of the paths $s_{1-}'$ or $s_{1+}'$ such that
both $u_1$ and $u_2$ are subpaths of~$z^{-1}$.
Then $\ell (z)$ is counter-regular.
Suppose $u_2$ precedes $u_1$ as a subpath of~$z^{-1}$.
Let $p$ be the path such that $u_1^{-1}pu_2^{-1}$ is a subpath of~$z$.
Then $\ell (u_1^{-1}pu_2^{-1})$ is counter-regular.
Therefore, $\ell (u_2^{-1}v^{-1}u_1^{-1})$ is regular.
The cycle $\langle pu_2^{-1}v^{-1}u_1^{-1}\rangle$ is the contour of
a disc subdiagram of $\Delta$ not containing the face $\Pi_1$
but containing the face~$\Pi_2$.
This contradicts the minimality of~$\Delta$.
Therefore, $u_1$ precedes $u_2$ as a subpath of~$z^{-1}$.
Let $p$ be the path such that $u_2^{-1}pu_1^{-1}$ is a subpath of~$z$.
Then $\ell (u_2^{-1}pu_1^{-1})$ is counter-regular.
Therefore, $\ell (u_1vu_2)$ is regular.
The cycle $\langle pv\rangle$ is the contour of
a disc subdiagram of $\Delta$ not containing the faces
$\Pi_1$ and~$\Pi_2$.
By Lemma~\ref{lemma:(5.4)}
(by the minimality of $\Delta$), this subdiagram is degenerate.
Since $z$ and $v$ are reduced, $v=p^{-1}$.
Therefore, $u_1vu_2$ is a subpath of~$z^{-1}$.
This gives a contradiction.

Case~3: $u_1$ is a subpath of $s_{1-}^{\prime-1}$, and
$u_2$ is a subpath of~$s_{1+}^{\prime-1}$.
Let $p$ be the path such that $u_1^{-1}p^{-1}u_2^{-1}$ is
a subpath of~$s_1$.
Then $\ell (u_2pu_1)$ is regular.
Therefore, $\ell (u_1vu_2)$ is counter-regular.
The cycle $\langle pu_1vu_2\rangle$ is the contour of a disc subdiagram
of $\Delta$ containing the face $\Pi_1$ but not containing
the face~$\Pi_2$.
This contradicts the minimality of $\Delta$ (see Lemma~\ref{lemma:(5.4)}).

\begin{figure}
\includegraphics{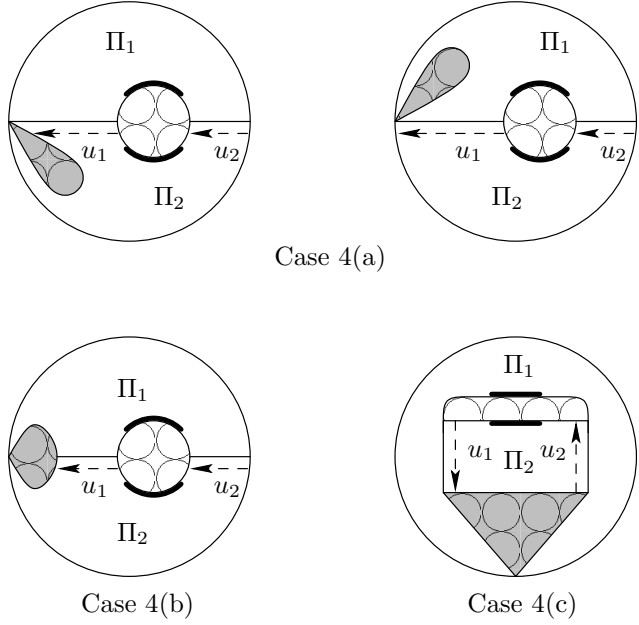}
\caption{Special subcases of Case 4, Lemma \ref{lemma:(5.10)}.}
\label{figure:05}
\end{figure}

Case~4: $u_1$ is a subpath of $s_{1+}^{\prime-1}$, and
$u_2$ is a subpath of~$s_{1-}^{\prime-1}$.
Since the labels of $s_{1-}'$ and $s_{1+}'$ have disjoint sets of
basic letters, $\ell (s_2)$ is regular.
Let $p$ be the path such that $u_2^{-1}p^{-1}u_1^{-1}$
is a subpath of~$s_1$.
Then $\ell (u_1vu_2)$ and $\ell (u_1pu_2)$ are regular.
If $v$ is trivial, then the cyclic path $pv^{-1}=p$ is cyclically reduced
because $\ell (p)$ is regular.
If $v$ is nontrivial, then the cyclic path $pv^{-1}$ is cyclically reduced
by the maximality of $u_1$ and $u_2$, and
because the pathes $v$ and $p$ are reduced.
Suppose the path $pv^{-1}$ is not simple.
Then some subpath of some cyclic shift of $vp^{-1}$ is a representative of
the contour of a proper simple disc subdiagram of~$\Delta$.
Let $\Delta'$ be such a subdiagram.
At least one of the following three subcases takes place:
\begin{enumerate}
\item[(a)]
	there is a subpath of $v$ or $p^{-1}$ which represents
	$\partial\Delta'$, or
\item[(b)]
	there are paths $v'$ and $p'$ such that
	$\langle v'p^{\prime-1}\rangle=\partial\Delta'$, and
	the paths $v'$ and $p'$ are either
	terminal subpaths of $v$ and $p$, respectively, or
	initial subpaths of $v$ and $p$, respectively, or
\item[(c)]
	there are an initial subpath $p_{-}'$ and a terminal subpath
	$p_{+}'$ of the path $p$ such that
$\langle vp_{+}^{\prime-1}p_{-}^{\prime-1}\rangle=\partial\Delta'$.
\end{enumerate}
In subcase (c) note that $\ell (p_{-}'p_{+}')$ is regular.
Lemma~\ref{lemma:(5.4)} easily yields a contradiction in each subcase.
Hence, the cyclic path $pv^{-1}$ is simple.
Therefore, it represents the contour of
a simple disc subdiagram of $\Delta$ containing
the faces $\Pi_1$ and~$\Pi_2$.
By Lemma~\ref{lemma:(5.4)}, $\langle pv^{-1}\rangle=\partial\Delta$.
Therefore, $v=q_1^{-1}$, and $u_1vu_2$ is a subpath of
$(s_{1-}'q_1s_{1+}')^{-1}$.
This gives a contradiction.

\begin{figure}
\includegraphics{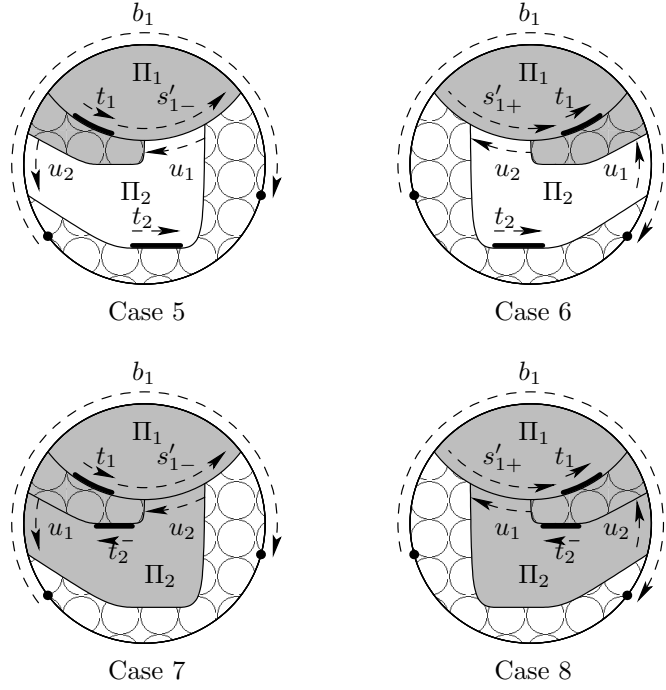}
\caption{Cases 5--8, Lemma \ref{lemma:(5.10)}.}
\label{figure:06}
\end{figure}

Case~5: $u_1$ is a subpath of $s_{1-}^{\prime-1}$, and
$u_2$ is a subpath of~$y_-^{-1}$.
Let $p$ be the path such that $u_2^{-1}pu_1$ is a subpath of
$y_-s_1^{\prime-1}s_{1-}^{\prime-1}$.
Then $\ell (u_2^{-1}pu_1)$ is regular.
Therefore, $\ell (u_1vu_2)$ is counter-regular.
The cycle $\langle pu_1v\rangle$ is the contour of
a disc subdiagram of $\Delta$ containing the face $\Pi_1$ but
not containing the face~$\Pi_2$.
This contradicts the minimality of $\Delta$ (see Lemma~\ref{lemma:(5.4)}).

Case~6: $u_1$ is a subpath of $y_+^{-1}$, and
$u_2$ is a subpath of~$s_{1+}^{\prime-1}$.
Let $p$ be the path such that $u_2pu_1^{-1}$ is a subpath of
$s_{1+}^{\prime-1}s_1^{\prime-1}y_+$.
Then $\ell (u_2pu_1^{-1})$ is regular.
Therefore, $\ell (u_1vu_2)$ is counter-regular.
The cycle $\langle pvu_2\rangle$ is the contour of
a disc subdiagram of $\Delta$ containing the face $\Pi_1$ but
not containing the face~$\Pi_2$.
This contradicts the minimality of~$\Delta$.

Case~7: $u_1$ is a subpath of $y_-^{-1}$, and
$u_2$ is a subpath of~$s_{1-}^{\prime-1}$.
Let $p$ be the path such that
$u_1^{-1}pu_2$ is a subpath of $y_-s_1^{\prime-1}s_{1-}^{\prime-1}$.
The cyclic path $pv^{-1}u_1^{-1}$ is cyclically reduced.
Suppose it is not simple.
Then at least one of the following two subcases takes place:
\begin{enumerate}
\item[(a)]
	there is a subpath of $v$ or $p^{-1}$ which represents the contour of
	a proper simple disc subdiagram of $\Delta$, or
\item[(b)]
	there are paths $v'$ and $p'$ such that
	$\langle v'p^{\prime-1}\rangle$ is the contour of a proper simple
	disc subdiagram of $\Delta$, and
	$v'$ and $p'$ are terminal subpaths of $v$ and $p$, respectively.
\end{enumerate}
Lemma~\ref{lemma:(5.4)} yields a contradiction in both subcases.
Hence, the cyclic path $pv^{-1}u_1^{-1}$ is simple.
Therefore, it represents the contour of
a disc subdiagram of $\Delta$ containing the faces $\Pi_1$ and~$\Pi_2$.
The label of $u_1^{-1}pu_2$ is regular.
Therefore, $\ell (u_2^{-1}v^{-1}u_1^{-1})$ is counter-regular.
By Lemma~\ref{lemma:(5.4)},
$\langle pv^{-1}u_1^{-1}\rangle=\partial\Delta$.
Therefore, $u_1vu_2$ is a subpath of $(s_{1-}'q_1)^{-1}$.
This gives a contradiction.

Case~8: $u_1$ is a subpath of $s_{1+}^{\prime-1}$, and
$u_2$ is a subpath of~$y_+^{-1}$.
Arguing as in Case~7, one obtains that $u_1vu_2$ is a subpath of
$(q_1s_{1+}')^{-1}$.
This gives a contradiction.

\begin{figure}
\includegraphics{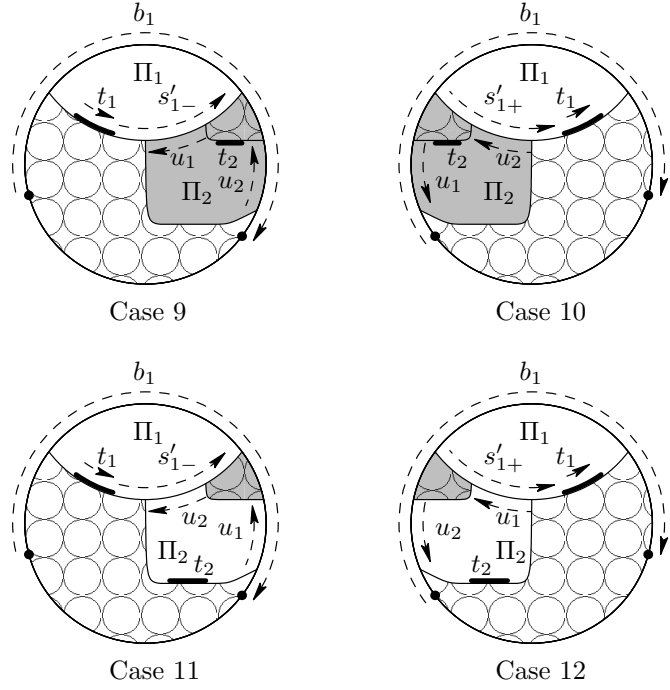}
\caption{Cases 9--12, Lemma \ref{lemma:(5.10)}.}
\label{figure:07}
\end{figure}

Case~9: $u_1$ is a subpath of $s_{1-}^{\prime-1}$, and
$u_2$ is a subpath of~$y_+^{-1}$.
Let $p_1$ and $p_2$ be the paths such that
$u_1^{-1}p_1$ is a terminal subpath of $s_{1-}'$,
and $p_2u_2^{-1}$ is an initial subpath of~$y_+$.
Then $\ell (u_1^{-1}p_1)$ is counter-regular,
$\ell (p_2u_2^{-1})$ is regular, and
$\ell (u_2^{-1}v^{-1}u_1^{-1})$ is regular or counter-regular.
If $\ell (u_2^{-1}v^{-1}u_1^{-1})$ is regular,
then $\ell (p_2u_2^{-1}v^{-1}u_1^{-1})$ is regular.
If $\ell (u_2^{-1}v^{-1}u_1^{-1})$ is counter-regular,
then $\ell (u_2^{-1}v^{-1}u_1^{-1}p_1)$ is counter-regular.
The cycle $\langle p_1p_2u_2^{-1}v^{-1}u_1^{-1}\rangle$ is the contour of
a disc subdiagram of $\Delta$ containing the face $\Pi_2$ but
not containing the face~$\Pi_1$.
This contradicts the minimality of~$\Delta$.

Case~10: $u_1$ is a subpath of $y_-^{-1}$, and
$u_2$ is a subpath of~$s_{1+}^{\prime-1}$.
Contradiction is obtained as in Case~9.

Case~11: $u_1$ is a subpath of $y_+^{-1}$, and
$u_2$ is a subpath of~$s_{1-}^{\prime-1}$.
Let $p_1$ and $p_2$ be the paths such that
$u_2^{-1}p_1$ is a terminal subpath of $s_{1-}'$, and
$p_2u_1^{-1}$ is an initial subpath of~$y_+$.
Then $\ell (u_2^{-1}p_1)$ is counter-regular,
$\ell (p_2u_1^{-1})$ is regular, and
$\ell (u_1vu_2)$ is regular or counter-regular.
If $\ell (u_1vu_2)$ is regular,
then the reduced form of $\ell (p_2v)$ is regular.
If $\ell (u_1vu_2)$ is counter-regular,
then the reduced form of $\ell (vp_1)$ is counter-regular.
The cycle $\langle p_1p_2v\rangle$ is the contour of
a disc subdiagram of $\Delta$ not containing
the faces $\Pi_1$ and~$\Pi_2$.
By Lemma~\ref{lemma:(5.4)}, $\langle p_1p_2v\rangle$ is the contour of
a degenerate disc subdiagram.
Since $s_{1-}'y_+$ and $v$ are reduced, $v=(p_1p_2)^{-1}$.
Therefore, $u_1vu_2$ is a subpath of~$(s_{1-}'y_+)^{-1}$.
This gives a contradiction.

Case~12: $u_1$ is a subpath of $s_{1+}^{\prime-1}$, and
$u_2$ is a subpath of~$y_-^{-1}$.
Arguing as in Case~11, one obtains that $u_1vu_2$ is a subpath
of~$(y_-s_{1+}')^{-1}$.
This gives a contradiction.

\begin{figure}
\includegraphics{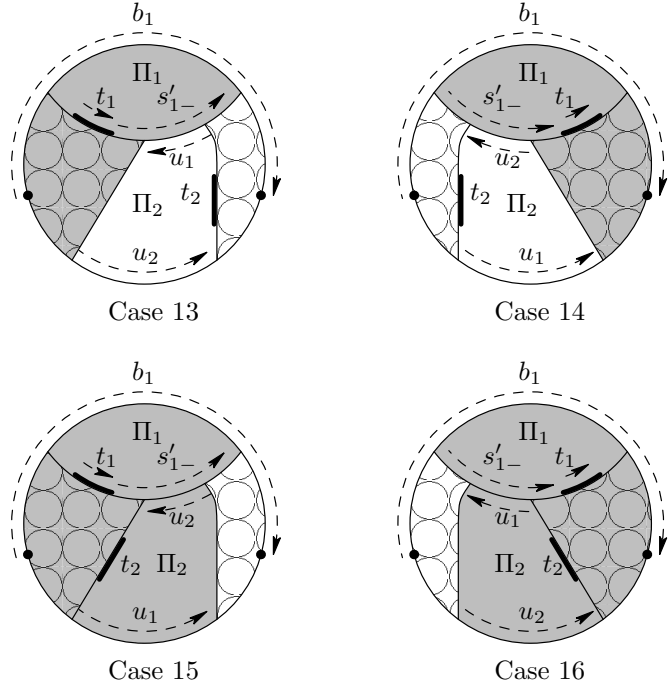}
\caption{Cases 13--16, Lemma \ref{lemma:(5.10)}.}
\label{figure:08}
\end{figure}

Case~13: $u_1$ is a subpath of $s_{1-}^{\prime-1}$, and
$u_2$ is a subpath of~$b_2^{\prime-1}$.
Let $p_1$ and $p_2$ be the paths such that
$p_1u_1$ is an initial subpath of $y_-s_1^{\prime-1}s_{1-}^{\prime-1}$,
and $u_2^{-1}p_2$ is a terminal subpath of~$b_2'$.
Then $\ell (p_1u_1)$ is regular,
$\ell (u_2^{-1}p_2)$ is counter-regular, and
$\ell (u_1vu_2)$ is regular or counter-regular.
If $\ell (u_1vu_2)$ is regular,
then $\ell (p_1u_1v)$ is regular.
If $\ell (u_1vu_2)$ is counter-regular,
then the reduced form of $\ell (u_1vp_2)$ is counter-regular.
The cycle $\langle p_2p_1u_1v\rangle$ is the contour of
a disc subdiagram of $\Delta$ containing the face $\Pi_1$ but
not containing the face~$\Pi_2$.
This contradicts the minimality of~$\Delta$.

Case~14: $u_1$ is a subpath of $b_2^{-1}$, and
$u_2$ is a subpath of~$s_{1+}^{\prime-1}$.
Contradiction is obtained as in Case~13.

Case~15: $u_1$ is a subpath of $b_2^{\prime-1}$, and
$u_2$ is a subpath of~$s_{1-}^{\prime-1}$.
Let $p_1$ and $p_2$ be the paths such that
$u_1^{-1}p_1$ is a terminal subpath of $b_2'$, and
$p_2u_2$ is an initial subpath of $y_-s_1^{\prime-1}s_{1-}^{\prime-1}$.
The cyclic path $p_1p_2v^{-1}u_1^{-1}$ is cyclically reduced.
Suppose it is not simple.
Then at least one of the following two subcases takes place:
\begin{enumerate}
\item[(a)]
	there is a subpath of $v$ or $p_1^{-1}$ which represents
	the contour of a proper simple disc subdiagram of $\Delta$, or
\item[(b)]
	there are paths $v'$ and $p_2'$ such that
	$\langle v'p_2^{\prime-1}\rangle$ is the contour of a proper simple
	disc subdiagram of $\Delta$, and
	$v'$ and $p_2'$ are terminal subpaths of $v$ and $p_2$, respectively.
\end{enumerate}
Lemma~\ref{lemma:(5.4)} yields a contradiction in both subcases.
Hence, $p_1p_2v^{-1}u_1^{-1}$ is a simple path.
Therefore, it represents the contour of
a disc subdiagram of $\Delta$ containing the faces
$\Pi_1$ and~$\Pi_2$.
The label of $u_1^{-1}p_1$ is counter-regular,
the label of $p_2u_2$ is regular, and
the label of $u_2^{-1}v^{-1}u_1^{-1}$ is regular or counter-regular.
If $\ell (u_2^{-1}v^{-1}u_1^{-1})$ is regular,
then the reduced form of $\ell (p_2v^{-1}u_1^{-1})$ is regular.
If $\ell (u_2^{-1}v^{-1}u_1^{-1})$ is counter-regular,
then $\ell (v^{-1}u_1^{-1}p_1)$ is counter-regular.
By Lemma~\ref{lemma:(5.4)},
$\langle p_1p_2v^{-1}u_1^{-1}\rangle=\partial\Delta$.
Therefore, $u_1vu_2$ is a subpath of $(s_{1-}'q_1)^{-1}$.
This gives a contradiction.

Case~16: $u_1$ is a subpath of $s_{1+}^{\prime-1}$, and
$u_2$ is a subpath of~$b_2^{\prime-1}$.
Arguing as in Case~15, one obtains that $u_1vu_2$ is a subpath of
$(q_1s_{1+}')^{-1}$.
This gives a contradiction.

A contradiction is obtained in each of the considered cases,
and no other case is possible.
Thus, $s_2$ and $(s_{1-}'q_1s_{1+}')^{-1}$ have exactly one
maximal common nontrivial subpath.
\end{proof}

\begin{lemma}
\label{lemma:(5.11)}
If\/ $\Pi$ is a face of\/ $\Delta$\textup,
the rank of every other face of\/ $\Delta$ is less than
the rank of\/ $\Pi$\textup,
$s$ is the maximal selected subpath
of\/ $\partial\Pi$\textup, and\/
$e_1$ and\/ $e_n$ are\textup, respectively\textup,
the initial and the terminal subpaths of\/ $s$
of length\/ $|s|/n$\textup,
then
all oriented edges of\/ $e_1$ or
all oriented edges of\/ $e_n$
are internal in\/~$\Delta$\textup.
\end{lemma}

\begin{proof}
Observe that proving that all oriented edges of $e_1$ or
all oriented edges of $e_n$ are internal in $\Delta$ is equivalent
to proving that all oriented edges of $e_n^{-1}$ or
all oriented edges of $e_1^{-1}$ are internal in
the mirror copy of~$\Delta$.

Suppose $\Pi$, $s$, $e_1$, $e_n$ are such as in the hypotheses
of the lemma.
Assume that $\ell (s)$ is regular
(if it is not, pass from $\Delta$, $s$, $e_1$, $e_n$
to the mirror copy of $\Delta$, $s^{-1}$, $e_n^{-1}$, $e_1^{-1}$,
respectively).

\begin{figure}
\includegraphics{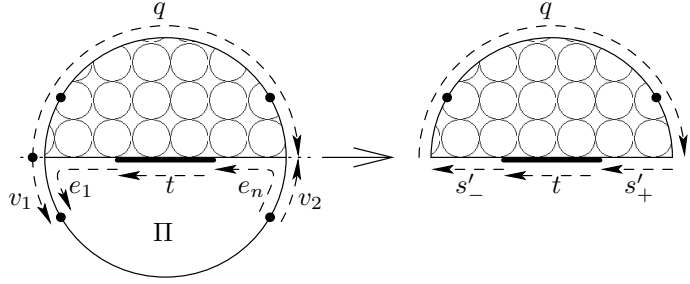}
\caption{The face $\Pi$ in $\Delta$ and the subdiagram $\Delta'$,
Lemma \ref{lemma:(5.11)}.}
\label{figure:09}
\end{figure}

Suppose that some oriented edge of $e_1$ and some oriented edge of $e_n$
are external in $\Delta$ (see Fig.~\ref{figure:09}).

Let $s'$ be the maximal selected external oriented arc of $\Delta$
that is a subpath of~$s$.
Let $s_{-}'$ and $s_{+}'$ be such paths that $s=s_{-}'s's_{+}'$.
Then $s_{-}'$ is a proper initial subpath of $e_1$, and
$s_{+}'$ is a proper terminal subpath of~$e_n$.
Let $v_1$ and $v_2$ be the initial and the terminal subpaths of $s'$
such that $s_{-}'v_1=e_1$ and $v_2s_{+}'=e_n$.
Let $t$ be the path such that $\langle st\rangle=\partial\Pi$.
Let $q$ be the path such that
$\langle qs^{\prime-1}\rangle=\partial\Delta$.
Let $\Delta'$ be the disc subdiagram of $\Delta$ obtained
by removing the face $\Pi$ and all edges that lie on $s'$
together with all intermediate vertices of~$s'$.
The contour of $\Delta'$ is $\langle s_{+}'ts_{-}'q\rangle$.

Let $j=\rank(\Pi)$.
Then $\ell (e_1)=x_1^{m_j}$ and $\ell (e_n)=x_n^{m_j}$.
The label of the first oriented edge of $s'$ is $x_1$;
the label of the last oriented edge of $s'$ is~$x_n$.
Therefore, $s'$ is a subpath of~$b_2^{-1}$.
Note that $v_1^{-1}$ is an initial subpath of $l_1$,
$v_2^{-1}$ is a terminal subpath of $l_2$,
and $v_1^{-1}qv_2^{-1}=b_1$.
The reduced form of $\ell (s_{-}'qs_{+}')$ is regular
since $\ell (q)$ is regular,
$\ell (s_{-}')$ is a power of $x_1$,
and $\ell (s_{+}')$ is a power of~$x_n$.
As it has been assumed, $\rank(\Pi')<j$ for every face $\Pi'$
of~$\Delta'$.
Hence, the group word $w_j=\ell (t^{-1})$
equals a regular word (the reduced form of $\ell (s_{-}'qs_{+}')$)
modulo the relations $r=1$, $r\in\mathcal R_{j-1}$.
This contradicts the choice of~$w_i$.
\end{proof}

\begin{lemma}
\label{lemma:(5.12)}
There are two distinct faces\/ $\Pi_1$ and\/ $\Pi_2$\textup, and
paths\/
$s_1$\textup, $s_1'$\textup, $s_{1-}'$\textup, $s_{1+}'$\textup,
$s_2$\textup, $s_2'$\textup, $s_{2-}'$\textup, $s_{2+}'$\textup,
$z_1$\textup, $z_2$ in\/ $\Delta$ such that\/\textup:
\begin{itemize}
\item
	$s_1$ and\/ $s_2$ are the maximal selected subpaths
	of\/ $\partial\Pi_1$ and\/ $\partial\Pi_2$\textup,
	respectively\/\textup;
\item
	$\ell (s_1)$ is counter-regular\textup,
	$\ell (s_2)$ is regular\/\textup;
\item
	$s_i=s_{i-}'s_i's_{i+}'$ for both\/ $i=1$ and\/ $i=2$\textup;
\item
	both\/ $s_1'$ and\/ $s_2'$ are maximal selected external
	oriented arcs of\/~$\Delta$\textup;
\item
	$\langle z_1s_1^{\prime-1}z_2s_2^{\prime-1}\rangle
	=\partial\Delta$\textup;
\item
	all basic letters of\/ $\ell (s_{2-}'z_1s_{1+}')$
	are in\/ $\{x_1,\dots,x_{21}\}$\textup, and
	all basic letters of\/ $\ell (s_{1-}'z_2s_{2+}')$
	are in\/ $\{x_{n-20},\dots,x_n\}$\textup.
\end{itemize}
\end{lemma}

\begin{proof}
Observe that this lemma is equivalent to the analogous lemma about
the mirror copy of $\Delta$,
\ie, to the statement obtained from this lemma by substituting
``the mirror copy of $\Delta$'' for ``$\Delta$.''

Let $\Pi_1$ be a face of $\Delta$ such as in Lemma~\ref{lemma:(5.9)}.
Let $s_1$ be the maximal selected subpath of $\partial\Pi_1$.
If $\ell (s_1)$ is regular, pass from $\Delta$ and $s_1$
to the mirror copy of $\Delta$ and $s_1^{-1}$, respectively.
Hence, $\ell (s_1)$ shall be assumed to be counter-regular.
Let $s_1'$ be the maximal selected external oriented arc of $\Delta$
that is a subpath of~$s_1$.
Let $s_{1-}'$ and $s_{1+}'$ be the paths such that
$s_1=s_{1-}'s_1's_{1+}'$.
Since all of the basic letters $x_{21}$, $x_{22}$, \dots, $x_{n-20}$
occur in $\ell (s_1')$, it follows that
all basic letters of $\ell (s_{1-}')$
are in $\{x_{n-20},\dots,x_n\}$, and
all basic letters of $\ell (s_{1+}')$
are in $\{x_1,\dots,x_{21}\}$.

Let $t_1$ be the path such that $\langle s_1t_1\rangle=\partial\Pi_1$.
Let $q_1$ be the path such that
$\langle s_1^{\prime-1}q_1\rangle=\partial\Delta$.
Let $\Delta'$ be the disc S-subdiagram of $\Delta$ obtained
by removing the face $\Pi_1$ and all edges that lie on $s_1'$
together with all intermediate vertices of~$s_1'$.
The contour of $\Delta'$ is $\langle s_{1+}'t_1s_{1-}'q_1\rangle$.

It shall be shown that there is a face $\Pi$ of $\Delta'$ such that
$\Pi$ is incident to no fewer than $(1-21/n)|\partial\Pi|$
selected external edges of $\Delta'$ that do not lie on~$t_1$.
Consider 2 cases:

Case~1: there is a face $\Pi$ in $\Delta'$ such that
$\rank(\Pi)\ge\rank(\Pi_1)$.
Then let $\Phi$ be the maximal simple disc S-subdiagram of $\Delta'$
that contains such a face~$\Pi$.
Let $\Sigma$ be the sum of the degrees of all the faces of~$\Phi$.
Then $\Sigma\ge |\partial\Pi|\ge |\partial\Pi_1|$.
By Lemma~\ref{lemma:(5.3)}, the number of selected external edges
of $\Phi$ is at least
$(1-2\mu)\Sigma$.
Therefore, the number of selected external edges of $\Phi$
that do not lie on $t_1$ is at least
$$
(1-2\mu)\Sigma-\lambda_1 |\partial\Pi_1|\ge(1-2\mu-\lambda_1)\Sigma
>\frac{n-21}{n}\Sigma.
$$

Case~2: the rank of every face of $\Delta'$ is less than
the rank of~$\Pi_1$.
Let $e_{11}$ and $e_{1n}$ be respectively the initial and terminal
subpaths of $s_1$ of length $|s_1|/n$.
Then $\ell (e_{11})=x_n^{-m_{j_1}}$, and $\ell (e_{1n})=x_1^{-m_{j_1}}$.
Since the rank of $\Pi_1$ is greater than the rank of every other face
of $\Delta$, all oriented edges of $e_{11}$ or all oriented edges
of $e_{1n}$ are internal in $\Delta$ by Lemma~\ref{lemma:(5.11)}.
Let $e$ denote $e_{11}$ in the first case or $e_{1n}$ in the second.

If an edge lies on $e$ and is not incident to any face of $\Delta'$,
then it also lies on~$t_1$.
Indeed, both (mutually inverse) oriented edges corresponding
to such an edge must occur in $\partial\Pi_1$;
one of them is an oriented edge of $e$, and consequently
the other cannot be an oriented edge of $s_1$
(two mutually inverse group letters cannot both occur in $\ell (s_1)$)
and has to be an oriented edge of~$t_1$.
Since
$$
|t_1|=|w_{j_1}|<m_{j_1}=|e|
$$
(see \thetag{\ref{display:(5.5)}}),
the ratio of the number of edges that lie on $t_1$
and are incident with faces of $\Delta'$
to the number of edges that lie on $e$ but are also incident with
faces of $\Delta'$ is not greater than
$$
\frac{|t_1|}{|e|}=\frac{|w_{j_1}|}{m_{j_1}}
\le\frac{n\lambda_1}{1-\lambda_1}
$$
(see \thetag{\ref{display:(5.4)}}).
Let $\Phi$ be a maximal simple disc S-subdiagram of $\Delta'$ such that
the ratio of the number of its (external) edges that lie on $t_1$
to the number of its (external) edges that lie on $e$
is not greater than
$$
\frac{n\lambda_1}{1-\lambda_1}.
$$
Let $\Sigma$ be the sum of the degrees of all the faces of~$\Phi$.
Since all edges that lie on $e$ are labelled with a same basic letter,
the number of edges of $\Phi$ that lie on $e$ is less than
$$
\Bigl(\lambda_1+\frac{1-\lambda_1}{n}\Bigr)\Sigma
$$
(because the number of edges with a same label incident
to a given face $\Pi$
of $\Delta$ is less than $(\lambda_1+(1-\lambda_1)/n)|\partial\Pi|$).
By Lemma~\ref{lemma:(5.3)}, the number of selected external
edges of $\Phi$ is at least
$(1-2\mu)\Sigma$.
Therefore, the number of selected external edges of $\Phi$
that do not lie on $t_1$ is greater than
$$
\Bigl(1-2\mu-\frac{n\lambda_1}{1-\lambda_1}
\Bigl(\lambda_1+\frac{1-\lambda_1}{n}\Bigr)\Bigr)\Sigma
=\Bigl(1-2\mu-\lambda_1-\frac{n\lambda_1^2}{1-\lambda_1}\Bigr)\Sigma
>\frac{n-21}{n}\Sigma
$$
(see \thetag{\ref{display:(5.2)}}).

In both cases, there is a face $\Pi$ in $\Phi$ such that $\Pi$
is incident to at least $(1-21/n)|\partial\Pi|$ selected external edges
of $\Phi$ (of $\Delta'$) that do not lie on~$t_1$.
Let $\Pi_2$ be such a face.
Let $s_2$ be the maximal selected subpath of $\partial\Pi_2$.
Then the paths $s_2$ and $(s_{1-}'q_1s_{1+}')^{-1}$ have at least
$(1-21/n)|\partial\Pi_2|$ common oriented edges.
Let $\tilde s_2$ be the maximal common nontrivial subpath of
$s_2$ and $(s_{1-}'q_1s_{1+}')^{-1}$
(see Lemma~\ref{lemma:(5.10)}).
Note that $\tilde s_2$ is a selected external oriented arc of~$\Delta'$.
Every selected external edge of $\Delta'$ incident to $\Pi_2$
lies either on $t_1$ or on~$\tilde s_2$.
Therefore,
$$
|\tilde s_2|\ge\frac{n-21}{n}|\partial\Pi_2|>\frac{n-21}{n}|s_2|.
$$
Therefore, $\ell (\tilde s_2)$ has at least $n-20$ distinct basic letters,
which implies that each of the basic letters $x_{21}$, \dots, $x_{n-20}$
occurs in it.
The basic letters $x_{22}$, \dots, $x_{n-21}$ do not occur on $s_{1-}'$
nor on $s_{1+}'$, but they occur in $\ell (\tilde s_2)$.
Hence, the paths $\tilde s_2$ and $q_1^{-1}$ have common oriented edges.
Therefore, $\Pi_2$ is incident to some selected external
edges of~$\Delta$.

Let $s_2'$ be the maximal selected external oriented arc of $\Delta$
that is a subpath of~$s_2$.
Then $s_2'$ is the maximal common subpath of $\tilde s_2$ and~$q_1^{-1}$.
All of the basic letters $x_{22}$, \dots, $x_{n-21}$
occur in $\ell (s_2')$
(since they occur in $\ell (\tilde s_2)$ but neither in $\ell (s_{1-}')$
nor in $\ell (s_{1+}')$).

Since $\ell (s_1)$ is counter-regular,
$s_1'$ is a subpath of~$b_1^{-1}$.
Therefore, $s_2'$ is a subpath of $b_2^{-1}$, and
$\ell (s_2)$ is regular.

Let $s_{2-}'$ and $s_{2+}'$ be the paths such that
$s_2=s_{2-}'s_2's_{2+}'$.
Let $z_1$ and $z_2$ be the paths such that
$\langle z_1s_1^{\prime-1}z_2s_2^{\prime-1}\rangle=\partial\Delta$.
It shall be proved that
all basic letters of $\ell (s_{2-}')$ and $\ell (z_1)$
are among $x_1$, \dots, $x_{21}$, and
all basic letters of $\ell (s_{2+}')$ and $\ell (z_2)$
are among $x_{n-20}$, \dots,~$x_n$.

First, consider $s_{2-}'$ and~$z_1$.

Case~1: $\tilde s_2$ has no common oriented edges
with~$s_{1+}^{\prime-1}$.
Then $s_2'$ is an initial subpath of~$\tilde s_2$.
Therefore, the basic letter $x_{21}$, along with
$x_{22}$, \dots, $x_{n-21}$,
occurs in $\ell (s_2')$.
Hence, all basic letters of $\ell (s_{2-}')$
are in $\{x_1,\dots,x_{21}\}$.
The label of every edge that lies on both $z_1$ and $b_1$ is
in $\{x_1,\dots,x_{21}\}$ since
the basic letter $x_{21}$ occurs in $\ell (s_1')$, and
the only maximal common subpath of $z_1s_1^{\prime-1}$ and $b_1$
is an initial subpath of~$b_1$.
The label of every edge that lies on both $z_1$ and $b_2$ is
in $\{x_1,\dots,x_{21}\}$ since
the basic letter $x_{21}$ occurs in $\ell (s_2')$, and
the only maximal common subpath of $s_2^{\prime-1}z_1$ and $b_2$
is a terminal subpath of~$b_2$.

Case~2: $\tilde s_2$ has at least one oriented edge
in common with~$s_{1+}^{\prime-1}$.
Then some nontrivial terminal subpath of $s_{1+}^{\prime-1}$ is
an initial subpath of $\tilde s_2$, and $z_1$ is trivial.
The terminal subpath of $s_{1+}^{\prime-1}$ that is an initial
subpath of $\tilde s_2$ is also a terminal subpath of~$s_{2-}'$.
%Since $\ell (s_{1+}^{\prime-1})$ and $\ell (s_{2-}')$ are both regular,
%and these paths have a common nontrivial terminal subpath,
%the set of basic letters that occur on $s_{2-}'$ is the same as
%the set of basic letters that occur on~$s_{1+}'$.
Since $s_{1+}^{\prime-1}$ and $s_{2-}'$ have a common nontrivial
terminal subpath, the sets of basic letters of their labels coincide.
Therefore, all basic letters of $\ell (s_{2-}')$
are in $\{x_1,\dots,x_{21}\}$.
The label of $z_1$ is empty.

Second, consider $s_{2+}'$ and $z_2$ in the same manner.

Case~1: $\tilde s_2$ has no common oriented edges
with~$s_{1-}^{\prime-1}$.
Then $s_2'$ is a terminal subpath of~$\tilde s_2$.
Therefore, the basic letter $x_{n-20}$,
along with $x_{22}$, \dots, $x_{n-21}$, occurs in $\ell (s_2')$.
Hence, all basic letters of $\ell (s_{2+}')$
are in $\{x_{n-20},\dots,x_n\}$.
The label of every edge that lies on both $z_2$ and $b_1$ is
in $\{x_{n-20},\dots,x_n\}$ since
the basic letter $x_{n-20}$ occurs in $\ell (s_1')$, and
the only maximal common subpath of $s_1^{\prime-1}z_2$ and $b_1$
is a terminal subpath of~$b_1$.
The label of every edge that lies on both $z_2$ and $b_2$ is
in $\{x_{n-20},\dots,x_n\}$ since
the basic letter $x_{n-20}$ occurs in $\ell (s_2')$, and
the only maximal common subpath of $z_2s_2^{\prime-1}$ and $b_2$
is an initial subpath of~$b_2$.

Case~2: $\tilde s_2$ has at least one oriented edge
in common with~$s_{1-}^{\prime-1}$.
Then some nontrivial initial subpath of $s_{1-}^{\prime-1}$ is
a terminal subpath of $\tilde s_2$, and $z_2$ is trivial.
The initial subpath of $s_{1-}^{\prime-1}$ that is a terminal
subpath of $\tilde s_2$ is also an initial subpath of~$s_{2+}'$.
%Since $\ell (s_{1-}^{\prime-1})$ and $\ell (s_{2+}')$ are both regular,
%and these paths have a common nontrivial initial subpath,
%the set of basic letters that occur on $s_{2+}'$ is the same as
%the set of basic letters that occur on~$s_{1-}'$.
Since $s_{1-}^{\prime-1}$ and $s_{2+}'$ have a common nontrivial
initial subpath, the sets of basic letters of their labels coincide.
Therefore, all basic letters of $\ell (s_{2+}')$
are in $\{x_{n-20},\dots,x_n\}$.
The label of $z_2$ is empty.

Clearly, the faces $\Pi_1$, $\Pi_2$ and the paths
$s_1$, $s_1'$, $s_{1-}'$, $s_{1+}'$,
$s_2$, $s_2'$, $s_{2-}'$, $s_{2+}'$, $z_1$, $z_2$
are desired ones.
\end{proof}

\begin{lemma}
\label{lemma:(5.13)}
The diagram\/ $\Delta$ has more than\/ $2$ faces\textup.
\end{lemma}

\begin{proof}
Note that $\Delta$ obviously has the same number of faces as
its mirror copy.

Suppose the statement is not true, \ie, $\Delta$ has no more
than $2$ faces.
(Then, by Lemma~\ref{lemma:(5.5)}, it has exactly $2$ faces.)

Let $\Pi_1$ and $\Pi_2$ be such faces of $\Delta$ as
in Lemma~\ref{lemma:(5.12)}.
Then the diagram $\Delta$ consists of the faces $\Pi_1$ and $\Pi_2$
attached to each other along a common arc.
If $\rank(\Pi_2)>\rank(\Pi_1)$, pass to the mirror copy of $\Delta$
and interchange the roles of $\Pi_1$ and~$\Pi_2$.
Hence, it shall be assumed that $\rank(\Pi_1)\ge\rank(\Pi_2)$.

Let $s_1$, $s_1'$, $s_{1-}'$, $s_{1+}'$,
$s_2$, $s_2'$, $s_{2-}'$, $s_{2+}'$ be the subpaths of
$\partial\Pi_1$ and $\partial\Pi_2$ such as in Lemma~\ref{lemma:(5.12)}.
Then all basic letters of the labels of
$s_{1-}'$, $s_{1+}'$, $s_{2-}'$, $s_{2+}'$
are in $\{x_1,\dots,x_{21}\}\cup\{x_{n-20},\dots,x_n\}$,
and $s_i'$ is a subpath of $b_i^{-1}$ for $i=1$ and $i=2$.

There are exactly two (mutually inverse) maximal internal oriented arcs
in~$\Delta$.
Denote the one that is a subpath of $\partial\Pi_1$ by~$u$.

Let $t_1$ and $t_2$ be the path such that
$\langle s_1t_1\rangle=\partial\Pi_1$ and
$\langle s_2t_2\rangle=\partial\Pi_2$.
Each of the faces $\Pi_1$ and $\Pi_2$ is incident to
exactly one maximal selected external arc of~$\Delta$.
Therefore, each of the paths $t_1$ and $t_2$ has a common vertex with~$u$.

Suppose $t_1$ and $t_2$ do not have common vertices.
Then either some (proper) nontrivial initial subpath of $s_1$ is
inverse to a (proper) initial subpath of $s_2$,
or some (proper) nontrivial terminal subpath of $s_1$ is
inverse to a (proper) terminal subpath of~$s_2$.
%Then either there exist two mutually inverse nontrivial paths which
%are proper initial subpaths of respectively $s_1$ and $s_2$,
%or there exist two mutually inverse nontrivial paths which
%are proper terminal subpaths of respectively $s_1$ and~$s_2$.
Both cases are impossible since
no proper nontrivial prefix of $\ell (s_2)$
can be a suffix of $\ell (s_1^{-1})$, and
no proper nontrivial suffix of $\ell (s_2)$
can be a prefix of $\ell (s_1^{-1})$.
(A contradiction may also be obtained by applying
Lemma~\ref{lemma:(5.10)}.)
Hence, the paths $t_1$ and $t_2$ have at least one common vertex.

Let $j_1=\rank(\Pi_1)$ and $j_2=\rank(\Pi_2)$.
Let $e_{i1}$ and $e_{in}$ be respectively the initial and the terminal
subpaths of $s_i$ of length $m_{j_i}=|s_i|/n$ for $i=1$ and $i=2$.
Then
$\ell (e_{11})=x_n^{-m_{j_1}}$,
$\ell (e_{1n})=x_1^{-m_{j_1}}$,
$\ell (e_{21})=x_1^{m_{j_2}}$, and
$\ell (e_{2n})=x_n^{m_{j_2}}$.
At most one of the two group letters $x_1^{\pm1}$ and at most one
of $x_n^{\pm1}$ can occur as labels of oriented edges of $\partial\Delta$.
Therefore, at least one of the paths $e_{1n}$ and $e_{21}$
and at least one of the paths $e_{11}$ and $e_{2n}$ do not have
external oriented edges in~$\Delta$.

Suppose $\rank(\Pi_1)=\rank(\Pi_2)$.
Then
$$
|e_{11}|=|e_{1n}|=|e_{21}|=|e_{2n}|=m_{j_1}\ge|w_{j_1}|=|t_1|=|t_2|.
$$
Recall that $w_{j_1}$ does not start with $x_1^{\pm1}$, does not end
with $x_n^{\pm1}$, and is not a letter power.
Therefore, every subpath of either $\partial\Pi_1$ or $\partial\Pi_2$
labelled with a letter power has length of at most~$m_{j_1}$.
Moreover, for every $i$,
there exist exactly one subpath of $\partial\Pi_1$
and exactly one subpath of $\partial\Pi_2$ of length $m_{j_1}$
labelled with powers of $x_i$
(\ie, with $x_i^{\pm m_{j_1}}$).
If $e_{1n}$ has no external oriented edges in $\Delta$, then
it is a subpath of~$u$.
If $e_{21}$ has no external oriented edges in $\Delta$, then
it is a subpath of~$u^{-1}$.
In either case $u$ has a subpath labelled with~$x_1^{-m_{j_1}}$,
and the pair $\{\Pi_1,\Pi_2\}$ is immediately cancellable.
This contradicts Lemma~\ref{lemma:(5.1)}.
Hence, $\rank(\Pi_1)>\rank(\Pi_2)$.

Let $\Delta'$ be the disc S-subdiagram of $\Delta$ obtained
by removing the face $\Pi_1$ and all edges and
intermediate vertices of~$s_1'$.
The only face of $\Delta'$ is~$\Pi_2$.

By Lemma~\ref{lemma:(5.11)}, all oriented edges of $e_{11}$ or
all oriented edges of $e_{1n}$ are internal in~$\Delta$.
Consider the following 5 cases (see Fig.~\ref{figure:10}):

(Notice that Case~3 is ``symmetric'' with Case~2,
and Case~5---with Case~4.)

\begin{figure}
\includegraphics{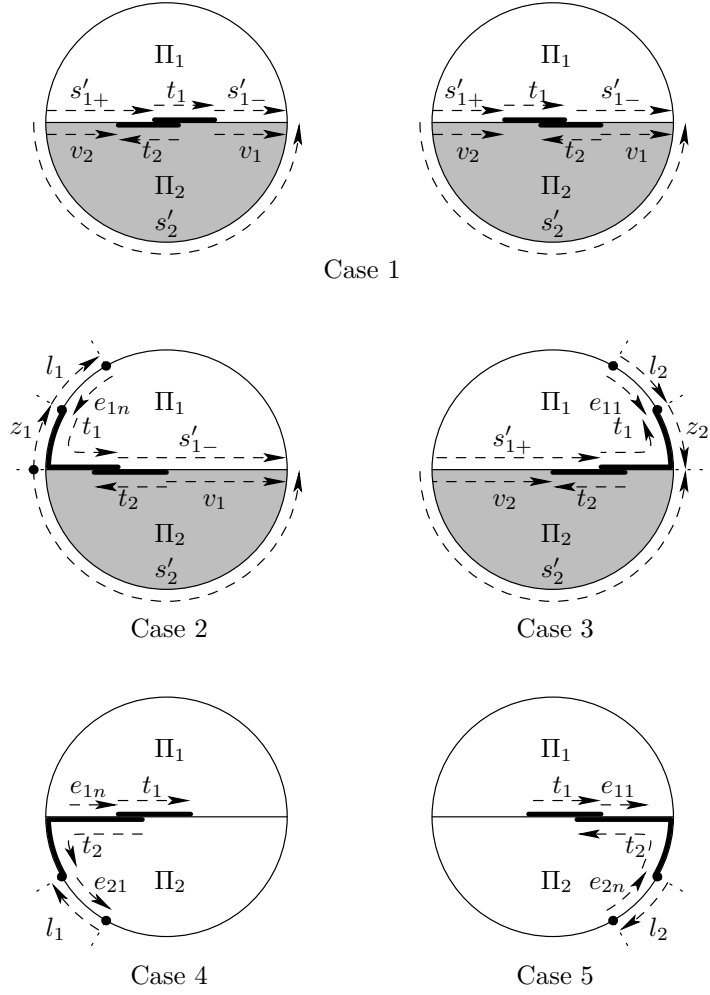}
\caption{Cases 1--5, Lemma \ref{lemma:(5.13)}.}
\label{figure:10}
\end{figure}

Case~1:
each of the paths $e_{11}$, $e_{1n}$, $e_{21}$, and $e_{2n}$ has
an internal oriented edge.
Then the paths $s_{1-}'$ and $s_{2+}^{\prime-1}$ have a common
nontrivial terminal subpath,
the paths $s_{1+}'$ and $s_{2-}^{\prime-1}$ have a common
nontrivial initial subpath,
and $\partial\Delta=\langle s_1^{\prime-1}s_2^{\prime-1}\rangle$.
Let $v_1$ be a common nontrivial terminal subpath of
$s_{1-}'$ and $s_{2+}^{\prime-1}$,
and let $v_2$ be a common nontrivial initial subpath of
$s_{1+}'$ and~$s_{2-}^{\prime-1}$.
The contour of $\Delta'$ is
$\langle s_{1+}'t_1s_{1-}'s_2^{\prime-1}\rangle$.
The label of $s_{1-}'s_2^{\prime-1}s_{1+}'$ is counter-regular
because the labels of $s_{1-}'$, $v_1s_2^{\prime-1}v_2$, and $s_{1+}'$
are such.
Hence, the group word $w_{j_1}=\ell (t_1)$
equals the regular word $\ell (s_{1-}'s_2^{\prime-1}s_{1+}')^{-1}$
modulo the relation $r_{j_2}=1$.
This contradicts the choice of~$w_{j_1}$.

Case~2:
all oriented edges of $e_{1n}$ are external
(consequently, all oriented edges of $e_{11}$ and $e_{21}$ are internal),
and some oriented edge of $e_{2n}$ is internal.
Then the paths $s_{1-}'$ and $s_{2+}^{\prime-1}$ have a common
nontrivial terminal subpath.
Let $v_1$ be a common nontrivial terminal subpath of
$s_{1-}'$ and~$s_{2+}^{\prime-1}$.
Let $z_1$ be the path such that
$\langle z_1s_1^{\prime-1}s_2^{\prime-1}\rangle=\partial\Delta$.
Then $z_1^{-1}$ is an initial subpath of~$t_1$.
Since $\ell (e_{1n})$ is a power of $x_1$, the path $e_{1n}^{-1}$
is a (terminal) subpath of~$b_2$.
Therefore, $s_2^{\prime-1}z_1e_{1n}^{-1}$ is a subpath of~$b_2$.
The contour of $\Delta'$ is $\langle z_1t_1s_{1-}'s_2^{\prime-1}\rangle$.
The label of $s_{1-}'s_2^{\prime-1}z_1$ is counter-regular
because the labels of $s_{1-}'$, $v_1s_2^{\prime-1}$, and
$s_2^{\prime-1}z_1$ are such.
Hence, the group word $w_{j_1}=\ell (t_1)$
equals the regular word $\ell (s_{1-}'s_2^{\prime-1}z_1)^{-1}$
modulo the relation $r_{j_2}=1$.
This contradicts the choice of~$w_{j_1}$.

Case~3:
all oriented edges of $e_{11}$ are external
(consequently, all oriented edges of $e_{1n}$ and $e_{2n}$ are internal),
and some oriented edge of $e_{21}$ is internal.
Then the paths $s_{1+}'$ and $s_{2-}^{\prime-1}$ have a common
nontrivial initial subpath.
Let $v_2$ be a common nontrivial initial subpath of
$s_{1+}'$ and~$s_{2-}^{\prime-1}$.
Let $z_2$ be the path such that
$\langle s_1^{\prime-1}z_2s_2^{\prime-1}\rangle=\partial\Delta$.
Then $z_2^{-1}$ is a terminal subpath of~$t_1$.
Since $\ell (e_{11})$ is a power of $x_n$, the path $e_{11}^{-1}$
is an (initial) subpath of~$b_2$.
Therefore, $e_{11}^{-1}z_2s_2^{\prime-1}$ is a subpath of~$b_2$.
The contour of $\Delta'$ is $\langle s_{1+}'t_1z_2s_2^{\prime-1}\rangle$.
The label of $z_2s_2^{\prime-1}s_{1+}'$ is counter-regular
because the labels of
$z_2s_2^{\prime-1}$, $s_2^{\prime-1}v_2$, and $s_{1+}'$ are such.
Hence, the group word $w_{j_1}=\ell (t_1)$
equals the regular word $\ell (z_2s_2^{\prime-1}s_{1+}')^{-1}$
modulo the relation $r_{j_2}=1$.
This contradicts the choice of~$w_{j_1}$.

Case~4:
all oriented edges of $e_{21}$ are external.
Consequently, all oriented edges of $e_{1n}$ are internal.
Then $e_{1n}$ is a terminal subpath of~$s_{1+}'$.
The path $s_{1+}^{\prime-1}$ is a subpath of $t_2$
(because $t_2$ has common vertices with both $e_{21}$ and~$t_1$).
Therefore, $|t_2|\ge|s_{1+}'|\ge|e_{1n}|=m_{j_1}>|w_{j_2}|=|t_2|$
(see \thetag{\ref{display:(5.5)}}).
This gives a contradiction.

Case~5:
all oriented edges of $e_{2n}$ are external.
Consequently, all oriented edges of $e_{11}$ are internal.
Then $e_{11}$ is an initial subpath of~$s_{1-}'$.
The path $s_{1-}^{\prime-1}$ is a subpath of $t_2$
(because $t_2$ has common vertices with both $e_{2n}$ and~$t_1$).
Therefore, $|t_2|\ge|s_{1-}'|\ge|e_{11}|=m_{j_1}>|w_{j_2}|=|t_2|$
(see \thetag{\ref{display:(5.5)}}).
This gives a contradiction.

No other case is possible.
Thus, $\Delta$ has more than $2$ faces.
\end{proof}

\begin{lemma}
\label{lemma:(5.14)}
Let\/ $\Pi_1$ and\/ $\Pi_2$ be distinct faces of\/ $\Delta$ such as in
Lemma\/~\textup{\ref{lemma:(5.12)}.}
Then the ranks of\/ $\Pi_1$ and\/ $\Pi_2$ are distinct\textup.
\end{lemma}

\begin{proof}
Let $\Pi_1$, $\Pi_2$, $s_1$, $s_1'$, $s_{1-}'$, $s_{1+}'$,
$s_2$, $s_2'$, $s_{2-}'$, $s_{2+}'$, $z_1$, $z_2$
be as in Lemma~\ref{lemma:(5.12)}.
Let $c_1=s_{2-}'z_1s_{1+}'$ and $c_2=s_{1-}'z_2s_{2+}'$.
Let $t_1$ and $t_2$ be the paths such that
$\langle s_1t_1\rangle=\partial\Pi_1$ and
$\langle s_2t_2\rangle=\partial\Pi_2$.
Let $\Delta'$ be the disc S-subdiagram of $\Delta$ obtained
by removing the faces $\Pi_1$ and $\Pi_2$ and
all edges and intermediate vertices of the paths $s_1'$ and $s_2'$
(see Fig.~\ref{figure:11}).

\begin{figure}
\includegraphics{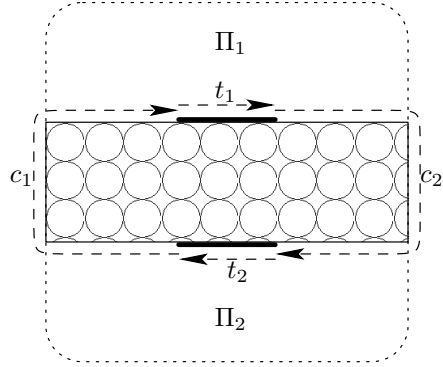}
\caption{The subdiagram $\Delta'$, Lemma \ref{lemma:(5.14)}.}
\label{figure:11}
\end{figure}

The contour of $\Delta'$ is $\langle c_1t_1c_2t_2\rangle$.
By Lemma~\ref{lemma:(5.13)}, the diagram $\Delta'$ is nondegenerate.
Therefore, by Lemma~\ref{lemma:(5.4)}, the label of
$c_1t_1c_2t_2$ is not freely trivial.

Suppose $\rank(\Pi_1)=\rank(\Pi_2)$.
Let $j=\rank(\Pi_1)=\rank(\Pi_2)$.
Recall that
$\ell (t_1)=w_j$, $\ell (t_2)=w_j^{-1}$,
all basic letters of $\ell (c_1)$ are in $\{x_1,\dots,x_{21}\}$, and
all basic letters of $\ell (c_2)$ are in $\{x_{n-20},\dots,x_n\}$.
Let $h_1=\ell (c_1^{-1})$ and $h_2=\ell (c_2)$.
Then $[h_1]=[w_j][h_2][w_j]^{-1}$.
Let $F$ be the subgroup of $G$ generated by
$\{a_1,\dots,a_{21}\}\cup\{a_{n-20},\dots,a_n\}$.
By Property~\ref{property:1}, the group $F$ is free of rank $42$
(since $42\le n-21$),
and elements of $F$ are conjugate in $G$ if and only if they are
conjugate in~$F$.
The elements $[h_1]$ and $[h_2]$ are in~$F$.
Moreover, $[h_1]$ belongs to the subgroup generated by
$\{a_1,\dots,a_{21}\}$,
and $[h_2]$ belongs to the subgroup generated by
$\{a_{n-20},\dots,a_n\}$.
Since $[h_1]$ and $[h_2]$ are conjugate in $G$, they are conjugate in $F$
and therefore $[h_1]=[h_2]=1$.
Hence, $\ell (c_1t_1c_2t_2)$ is freely trivial.
This gives a contradiction.
\end{proof}

\begin{lemma}
\label{lemma:(5.15)}
Let\/ $\Pi_1$ and\/ $\Pi_2$ be distinct faces of\/ $\Delta$\textup,
both incident to selected external edges of\/~$\Delta$\textup.
Suppose the rank of\/ $\Pi_1$ is greater than the rank
of every other face of\/~$\Delta$\textup.
Let\/ $s_1$ and\/ $s_2$ be the maximal selected subpaths
of\/ $\partial\Pi_1$ and\/ $\partial\Pi_2$\textup, respectively\textup.
Suppose the label of one of the paths\/ $s_1$ or\/ $s_2$ is
regular\textup, and the label of the other is counter-regular\textup.
Let\/ $s_1'$ and\/ $s_2'$ be the maximal selected external
oriented arcs of\/ $\Delta$ that are subpaths of\/
$s_1$ and\/ $s_2$\textup, respectively\textup.
Suppose there are at least\/ $2$ distinct basic letters in the label of
each of the paths\/ $s_1'$ and\/~$s_2'$\textup.
Let\/ $e_{11}$ and\/ $e_{1n}$ be respectively
the initial and the terminal subpaths of\/ $s_1$
of length\/~$|s_1|/n$\textup.
Then at least one of the paths\/ $e_{11}$ or\/ $e_{1n}$
has the property that every edge that lies on it
is internal in\/ $\Delta$ and does not lie on\/~$s_2$\textup.
\end{lemma}

\begin{proof}
Observe that it suffices to prove that
at least one of the paths $e_{1n}^{-1}$ or $e_{11}^{-1}$
has the property that every edge that lies on it
is internal in the mirror copy of $\Delta$ and does not lie on~$s_2^{-1}$.

Assume that $\ell (s_1)$ is counter-regular
(if it is not, pass from $\Delta$, $b_1$, $b_2$, $l_1$, $l_2$,
$s_1$, $s_2$, $s_1'$, $s_2'$, $e_{11}$, $e_{1n}$
to the mirror copy of
$\Delta$, $b_2^{-1}$, $b_1^{-1}$, $l_1^{-1}$, $l_2^{-1}$,
$s_1^{-1}$, $s_2^{-1}$, $s_1^{\prime-1}$, $s_2^{\prime-1}$,
$e_{1n}^{-1}$, $e_{11}^{-1}$, respectively).
Then $s_1'$ is a subpath of $b_1^{-1}$,
and $s_2'$ is a subpath of~$b_2^{-1}$.

Let $s_{1-}'$ and $s_{1+}'$ be the paths such that
$s_1=s_{1-}'s_1's_{1+}'$.
Let $t_1$ and $t_2$ be the paths such that
$\langle s_1t_1\rangle=\partial\Pi_1$ and
$\langle s_2t_2\rangle=\partial\Pi_2$.
Let $q_1$ be the path such that
$\langle s_1^{\prime-1}q_1\rangle=\partial\Delta$.
Let $\Delta'$ be the disc subdiagram of $\Delta$ obtained
by removing the face $\Pi_1$ and all edges that lie on $s_1'$
together with all intermediate vertices of~$s_1'$.
The contour of $\Delta'$ is $\langle s_{1+}'t_1s_{1-}'q_1\rangle$.
Let $e_{21}$ and $e_{2n}$ be respectively the initial and
terminal subpaths of $s_2$ of length $|s_2|/n$.

Denote the ranks of $\Pi_1$ and $\Pi_2$ by $j_1$ and $j_2$, respectively.
Then
\begin{gather*}
\alignedat 2
\ell (s_1)&=x_n^{-m_{j_1}}x_{n-1}^{-m_{j_1}}\dots x_1^{-m_{j_1}},&\qquad
\ell (t_1)&=w_{j_1},\\
\ell (s_2)&=x_1^{m_{j_2}}x_2^{m_{j_2}}\dots x_n^{m_{j_2}},&\qquad
\ell (t_2)&=w_{j_2}^{-1},
\endalignedat\\\vspace{1\jot}
\ell (e_{11})=x_n^{-m_{j_1}},\qquad
\ell (e_{1n})=x_1^{-m_{j_1}},\qquad
\ell (e_{21})=x_1^{m_{j_2}},\qquad
\ell (e_{2n})=x_n^{m_{j_2}}.
\end{gather*}

At most one of the two group letters $x_1^{\pm1}$ and at most one
of $x_n^{\pm1}$ can occur as labels of oriented edges of $\partial\Delta$
(the labels of the representatives of $\partial\Delta$
are cyclically reduced).
Therefore, at least one of the paths $e_{1n}$ or $e_{21}$
and at least one of the paths $e_{11}$ or $e_{2n}$ do not have
external oriented edges in~$\Delta$.
By Lemma~\ref{lemma:(5.11)}, all oriented edges of $e_{11}$ or
all oriented edges of $e_{1n}$ are internal in~$\Delta$.

Suppose the conclusion of the lemma does not hold.
Then some oriented edge of $e_{11}$
is either external in $\Delta$ or inverse to an oriented edge of $e_{2n}$,
and some oriented edge of $e_{1n}$
is either external in $\Delta$ or inverse to an oriented edge of~$e_{21}$.
Therefore, as follows from Lemma~\ref{lemma:(5.11)}, at least one of
the following 3 cases takes place:
\begin{itemize}
\item[$\circ$]
	some oriented edge of $e_{11}$
	is inverse to an oriented edge of $e_{2n}$, and
	some oriented edge of $e_{1n}$
	is inverse to an oriented edge of $e_{21}$, or
\item[$\circ$]
	some oriented edge of $e_{11}$
	is external in $\Delta$, and
	some oriented edge of $e_{1n}$
	is inverse to an oriented edge of $e_{21}$, or
\item[$\circ$]
	some oriented edge of $e_{11}$
	is inverse to an oriented edge of $e_{2n}$, and
	some oriented edge of $e_{1n}$
	is external in~$\Delta$.
\end{itemize}
In every case the path $s_2$ has a common oriented edge with
$s_{1-}^{\prime-1}$ or~$s_{1+}^{\prime-1}$.
Let $\tilde s_2$ be the maximal common nontrivial subpath
of $s_2$ and $(s_{1-}'q_1s_{1+}')^{-1}$
(see Lemma~\ref{lemma:(5.10)}).
Now, consider each of the 3 cases separately (see Fig.~\ref{figure:12}).

(Notice that Case~3 is ``symmetric'' with Case~2.)

\begin{figure}
\includegraphics{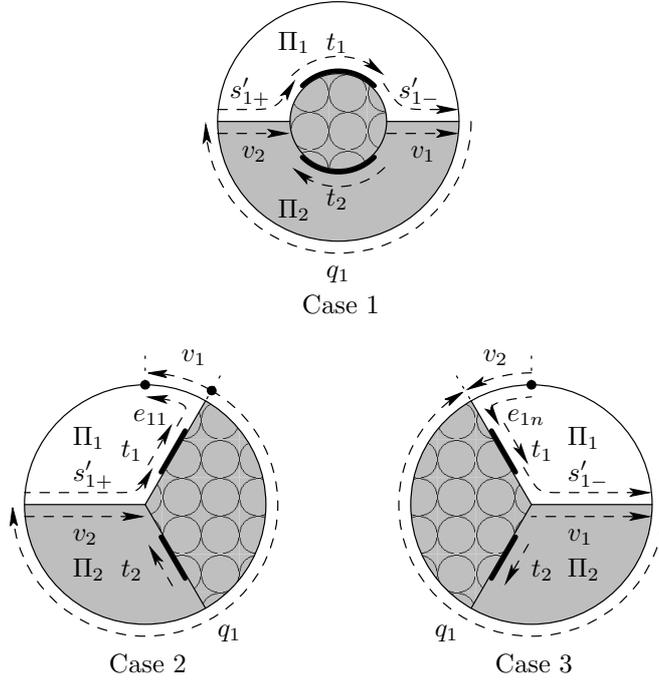}
\caption{Cases 1--3, Lemma \ref{lemma:(5.15)}.}
\label{figure:12}
\end{figure}

Case~1:
some oriented edge of $e_{11}$
is inverse to an oriented edge of $e_{2n}$, and
some oriented edge of $e_{1n}$
is inverse to an oriented edge of~$e_{21}$.
Then some nontrivial initial subpath of $\tilde s_2^{-1}$ is
a terminal subpath of $s_{1-}'$, and
some nontrivial terminal subpath of $\tilde s_2^{-1}$ is
an initial subpath of $s_{1+}'$
(because $\tilde s_2^{-1}$ is a subpath of $s_{1-}'q_1s_{1+}'$).
Let $v_1$ and $v_2$ be the paths such that $v_1q_1v_2=\tilde s_2^{-1}$.
Then $v_1$ is a nontrivial terminal subpath of $s_{1-}'$,
and $v_2$ is a nontrivial initial subpath of~$s_{1+}'$.
The label of $s_{1-}'q_1s_{1+}'$ is counter-regular
since the labels of $s_{1-}'$, $v_1q_1v_2$, and $s_{1+}'$ are such.
As it is assumed, $\rank(\Pi)<j_1$ for every face $\Pi$ of~$\Delta'$.
Recall that $\langle s_{1+}'t_1s_{1-}'q_1\rangle=\partial\Delta'$.
Hence, the group word $w_{j_1}=\ell (t_1)$
equals the regular word $\ell (s_{1-}'q_1s_{1+}')^{-1}$
modulo the relations $r=1$, $r\in\mathcal R_{j_1-1}$.
This contradicts the choice of~$w_{j_1}$.

Case~2:
some oriented edge of $e_{11}$ is external in $\Delta$, and
some oriented edge of $e_{1n}$ is inverse to an oriented edge of~$e_{21}$.
Then $s_{1-}'$ is a proper initial subpath of~$e_{11}$.
Let $v_1$ be the initial subpath of $s_1'$ such that $s_{1-}'v_1=e_{11}$.
Then $v_1^{-1}$ is a subpath of $l_2$ since $\ell (v_1^{-1})$ is
a nontrivial power of~$x_n$.
Some nontrivial terminal subpath of $\tilde s_2^{-1}$ is
an initial subpath of~$s_{1+}'$.
Therefore, $s_2^{\prime-1}$ is a terminal subpath of $q_1$, and
$v_1^{-1}q_1$ is a subpath of $b_2$
(because $s_2^{\prime-1}$ is a common subpath of $q_1$ and~$b_2$).
Let $v_2$ be the path such that $s_2^{\prime-1}v_2$ is
a terminal subpath of~$\tilde s_2^{-1}$.
Then $v_2$ is a nontrivial initial subpath of~$s_{1+}'$.
The label of $q_1s_{1+}'$ is counter-regular
since the labels of $q_1$, $s_2^{\prime-1}v_2$, and $s_{1+}'$ are such.
The reduced form of $\ell (s_{1-}'q_1s_{1+}')$ is counter-regular
since $\ell (s_{1-}')$ is a power of $x_n$, and
$\ell (q_1s_{1+}')$ is counter-regular.
As it is assumed, $\rank(\Pi)<j_1$ for every face $\Pi$ of~$\Delta'$.
Hence, the group word $w_{j_1}=\ell (t_1)$
equals a regular word (the reduced form of
$\ell (s_{1-}'q_1s_{1+}')^{-1}$)
modulo the relations $r=1$, $r\in\mathcal R_{j_1-1}$.
This contradicts the choice of~$w_{j_1}$.

Case~3:
some oriented edge of $e_{11}$ is inverse to an oriented edge of $e_{2n}$,
and some oriented edge of $e_{1n}$ is external in~$\Delta$.
Contradiction is obtained as in Case~2.

No other case is possible.
Thus, at least one of the paths $e_{11}$ or $e_{1n}$
has no common oriented edges with $s_2^{-1}$ and no external
oriented edges.
\end{proof}

\begin{lemma}
\label{lemma:(5.16)}
Let\/ $\Pi_1$\textup, $\Pi_2$\textup,
$s_1$\textup, $s_1'$\textup, $s_{1-}'$\textup, $s_{1+}'$\textup,
$s_2$\textup, $s_2'$\textup, $s_{2-}'$\textup, $s_{2+}'$\textup,
$z_1$\textup, $z_2$
be such as in Lemma\/~\textup{\ref{lemma:(5.12)}.}
Let\/ $c_1=s_{2-}'z_1s_{1+}'$ and\/ $c_2=s_{1-}'z_2s_{2+}'$\textup.
Then\/ $\Delta$ has a simple disc S-subdiagram\/ $\Phi$ such that\/
$\Pi_1$ and\/ $\Pi_2$ are not in\/ $\Phi$\textup, and the total
number of selected external edges of\/ $\Phi$
that lie on\/ $c_1$ or\/ $c_2$ is greater than
$$
\frac{n-21}{n}\!\sum_{\Pi\in\Phi(2)}\!|\partial\Pi|.
$$
\end{lemma}

\begin{proof}
Observe that to prove this lemma,
it suffices to prove the analogous lemma for the mirror copy of~$\Delta$.

By Lemma~\ref{lemma:(5.14)}, $\rank(\Pi_1)\ne\rank(\Pi_2)$.
If $\rank(\Pi_1)<\rank(\Pi_2)$, pass from $\Delta$,
$\Pi_1$, $\Pi_2$,
$s_1$, $s_{1-}'$, $s_1'$, $s_{1+}'$, $s_2$, $s_{2-}'$, $s_2'$, $s_{2+}'$,
$z_1$, $z_2$, $c_1$, $c_2$ to the mirror copy of $\Delta$,
$\Pi_2$, $\Pi_1$,
$s_2^{-1}$, $s_{2+}^{\prime-1}$, $s_2^{\prime-1}$, $s_{2-}^{\prime-1}$,
$s_1^{-1}$, $s_{1+}^{\prime-1}$, $s_1^{\prime-1}$, $s_{1-}^{\prime-1}$,
$z_1^{-1}$, $z_2^{-1}$, $c_1^{-1}$, $c_2^{-1}$, respectively.
Hence, assume that $j_1=\rank(\Pi_1)>\rank(\Pi_2)=j_2$.

Let $t_1$ and $t_2$ be the paths such that
$\langle s_1t_1\rangle=\partial\Pi_1$ and
$\langle s_2t_2\rangle=\partial\Pi_2$.
Let $\Delta'$ be the disc S-subdiagram of $\Delta$ obtained
by removing the faces $\Pi_1$ and $\Pi_2$ and all edges that lie
on $s_1'$ or $s_2'$
together with all intermediate vertices of $s_1'$ and~$s_2'$.
The contour of $\Delta'$ is $\langle c_1t_1c_2t_2\rangle$.
By Lemma~\ref{lemma:(5.13)}, the diagram $\Delta'$ is nondegenerate.

Consider 2 cases:

Case~1: there is a face $\Pi$ in $\Delta'$ such that
$\rank(\Pi)\ge\rank(\Pi_1)$.
Then let $\Phi$ be the maximal simple disc S-subdiagram of $\Delta'$
that contains such a face~$\Pi$.
Let $\Sigma$ be the sum of the degrees of all the faces of~$\Phi$.
Then $\Sigma\ge |\partial\Pi|\ge |\partial\Pi_1|>|\partial\Pi_2|$.
By Lemma~\ref{lemma:(5.3)}, the number of selected external
edges of $\Phi$ is at least $(1-2\mu)\Sigma$.
Therefore, the number of selected external edges of $\Phi$
that do not lie on $t_1$ nor on $t_2$ is at least
$$
(1-2\mu)\Sigma-\lambda_1 |\partial\Pi_1|-\lambda_1 |\partial\Pi_2|
>(1-2\mu-2\lambda_1)\Sigma>\frac{n-21}{n}\Sigma.
$$

Case~2: the rank of every face of $\Delta'$ is less than
the rank of~$\Pi_1$.
Let $e_{11}$ and $e_{1n}$ be respectively the initial and
the terminal subpaths of $s_1$ of length $m_{j_1}=|s_1|/n$.
By Lemma~\ref{lemma:(5.15)}, at least one of the paths
$e_{11}$ or $e_{1n}$ has the property that every edge that lies on
it either is incident to a face of $\Delta'$ or lies on $t_1$ or~$t_2$.
Let $e$ be one of the paths $e_{11}$ or $e_{1n}$ with this property.

Since
$$
|t_1|+|t_2|=|w_{j_1}|+|w_{j_2}|<\frac{1}{2}m_{j_1}+\frac{1}{2}m_{j_1}
=m_{j_1}=|e|
$$
(see \thetag{\ref{display:(5.5)}}),
the ratio of the number of edges that lie on $t_1$ or $t_2$
and are incident with faces of $\Delta'$
(consequently, do not lie on $e$)
to the number of edges that lie on $e$ and are incident with
faces of $\Delta'$ (equivalently, lie neither on $t_1$ nor on $t_2$)
is at most
$$
\frac{|t_1|+|t_2|}{|e|}=\frac{|w_{j_1}|+|w_{j_2}|}{m_{j_1}}
\le\frac{2n\lambda_1}{1-\lambda_1}
$$
(see \thetag{\ref{display:(5.4)}}).
Let $\Phi$ be a maximal simple disc S-subdiagram of $\Delta'$ such that
the ratio of the number of its (external) edges that lie on
$t_1$ or $t_2$ to the number of its (external) edges that lie on $e$
is at most
$$
\frac{2n\lambda_1}{1-\lambda_1}.
$$
Let $\Sigma$ be the sum of the degrees of all the faces of~$\Phi$.
Since all the edges that lie on $e$ are labelled with a same basic letter,
the number of edges of $\Phi$ that lie on $e$ is less than
$$
\Bigl(\lambda_1+\frac{1-\lambda_1}{n}\Bigr)\Sigma.
$$
By Lemma~\ref{lemma:(5.3)}, the number of selected external edges
of $\Phi$ is at least $(1-2\mu)\Sigma$.
Therefore, the number of selected external edges of $\Phi$
that neither lie on $t_1$ nor on $t_2$ is greater than
$$
\Bigl(1-2\mu-\frac{2n\lambda_1}{1-\lambda_1}
\Bigl(\lambda_1+\frac{1-\lambda_1}{n}\Bigr)\Bigr)\Sigma
=\Bigl(1-2\mu-2\lambda_1-\frac{2n\lambda_1^2}{1-\lambda_1}\Bigr)\Sigma
\ge\frac{n-21}{n}\Sigma
$$
(see \thetag{\ref{display:(5.2)}}).

Every external edge of $\Phi$ that neither lies on $t_1$ nor on $t_2$
lies on $c_1$ or~$c_2$.
Thus, the desired inequality holds in the both cases.
\end{proof}

%%%%%%%%%%%%%%%%%%%%%%%%%%%%%%%%%%%%%%%%%%%%%%%%%%%%%%%%%%%%%%%%%%%%%%%%%%%

\begin{proof}[Proof of Property\/ \textup{\ref{property:2}}]
The existence is obvious from the choice of relators.

Suppose that presentation of some element of $G$ in the form
$a_1^{k_1}\dots a_n^{k_n}$ is not unique.
Then there exist two distinct regular group words over $\mathfrak A$
whose values in $G$ coincide.
Consider such a pair of group words with the minimal sum of their lengths.
Consider a minimal deduction diagram for the equality of these group words
over $\langle\,\mathfrak A\,\|\,\mathcal R\,\rangle$.
More precisely, let $\Delta$ be a special disc S-diagram over
$\langle\,\mathfrak A\,\|\,\mathcal R\,\rangle$
whose contour has a representative of the form $p_1p_2$ such that
$\ell (p_1)$ and $\ell (p_2^{-1})$ are distinct regular group words,
and $\ell (p_1p_2)$ is cyclically reduced;
moreover, let $\Delta$ be such an S-diagram with the minimal possible
number of faces.

Choose faces $\Pi_1$ and $\Pi_2$ and paths $s_1$, $s_1'$, $s_{1-}'$, $s_{1+}'$,
$s_2$, $s_2'$, $s_{2-}'$, $s_{2+}'$, $z_1$, $z_2$ as
in Lemma~\ref{lemma:(5.12)}.
Denote the paths $s_{2-}'z_1s_{1+}'$ and $s_{1-}'z_2s_{2+}'$
by $c_1$ and $c_2$, respectively
(as in the hypotheses of Lemma~\ref{lemma:(5.16)}).
Apply Lemma \ref{lemma:(5.16)}:
let $\Phi$ be a simple disc S-subdiagram of $\Delta$ such that
the number of selected external edges of $\Phi$ that lie on $c_1$ or $c_2$
is greater than
$$
\frac{n-21}{n}\!\sum_{\Pi\in\Phi(2)}\!|\partial\Pi|.
$$
By the choice of $\Pi_1$ and $\Pi_2$,
all basic letters of $\ell (c_1)$ are in $\{x_1,\dots,x_{21}\}$, and
all basic letters of $\ell (c_2)$ are in $\{x_{n-20},\dots,x_n\}$.
Therefore, by Proposition~\ref{proposition:(5.2)},
the number of selected external edges of $\Phi$ that lie on $c_1$ or $c_2$
is less than
$$
\frac{42}{n}\sum_{\Pi\in\Phi(2)}\!|\partial\Pi|.
$$

Since $n-21\ge42$, one has a contradiction.
The uniqueness is proved.
\end{proof}

\begin{remark}
\label{remark:(5.1)}
The ``symmetry'' mentioned in the proofs of Lemmas
\ref{lemma:(5.8)}, \ref{lemma:(5.10)}, \ref{lemma:(5.13)},
and \ref{lemma:(5.15)}
is the symmetry between the S-diagram $\Delta$ and the S-diagram obtained
from the mirror copy of $\Delta$ by re-labelling according to
the following rule:
if the label of an oriented edge edge $e$ is $x_i^\sigma$,
then re-label $e$ with~$x_{n+1-i}^{-\sigma}$
(this gives a diagram not over the same presentation).
\end{remark}

The proof of the next property uses a homological argument.

\begin{property}
\label{property:3}
The group\/ $G$ is torsion-free\textup.
\end{property}

\begin{proof}
The group presentation $\langle\,\mathfrak A\,\|\,\mathcal R\,\rangle$ is
strongly aspherical by Corollary~\ref{corollary.proposition:2.(5.1)}
of Proposition~\ref{proposition:(5.1)};
therefore, it is aspherical in the sense
of~\cite{Olshanskii:1991:gdrg-eng}.
Since no element of $\mathcal R$ represents a proper power in the
free group on $\mathfrak A$,
the relation module $M$ of $\langle\,\mathfrak A\,\|\,\mathcal R\,\rangle$
is a free $G$-module by Corollary~32.1 in~\cite{Olshanskii:1991:gdrg-eng}.
Therefore, there exists a finite-length free resolution of $\mathbb Z$
over $\mathbb ZG$:
\begin{equation}
\label{display:(5.6)}
0\to M\to\bigoplus_{x\in\mathfrak A}\mathbb ZG
\to\mathbb ZG\to\mathbb Z\to0,
\end{equation}
where $\mathbb ZG$ and $\bigoplus_{x\in\mathfrak A}\mathbb ZG$ are
identified with
the $G$-modules of, respectively, $0$- and $1$-dimensional
cellular chains
of the Cayley complex of $\langle\,\mathfrak A\,\|\,\mathcal R\,\rangle$
(see~\cite{Brown:1994:cg}).

Now suppose that $G$ has torsion.
Let $H$ be a non-trivial finite cyclic subgroup of~$G$.
Every free $G$-module may be naturally regarded as a free $H$-module.
Hence, the resolution \thetag{\ref{display:(5.6)}} may be viewed
as a free resolution of $\mathbb Z$ over~$\mathbb ZH$.
This contradicts the fact that all odd-dimensional homology groups
of $H$ are nontrivial.
\end{proof}

For more results on group presentations with various
forms of asphericity,
and torsion in such groups, see \cite{ChiColHue:1981:agp},
\cite{Huebschmann:1979:ctagscg}, \cite{Huebschmann:1980:htcap},
and Theorems 32.1, 32.2 in \cite{Olshanskii:1991:gdrg-eng}.

Properties \ref{property:1}, \ref{property:2}, and \ref{property:3} of $G$
demonstrate that the answer to Bludov's question is negative.

\begin{remark}
\label{remark:(5.2)}
If $n$ was chosen to be less than $63$ but greater than $26$ ($26<n<63$),
then Properties \ref{property:1} and \ref{property:3} would still hold,
but the proof of Property~\ref{property:2} would not work.
\end{remark}

For every $i\in\mathbb N\cup\{0\}$,
let $G_i$ be the group defined by the presentation
$\langle\,\mathfrak A\,\|\,\mathcal R_i\,\rangle$
(see Section~\ref{section:group_construction}).
Also for every $i\in\mathbb N\cup\{0\}$,
let $\phi_i$ be the natural epimorphism $G_i\to G_{i+1}$.

\begin{property}
\label{property:4}
All groups\/ $G_0$\textup, $G_1$\textup, $G_2$\textup, \dots
are hyperbolic\textup, and\/
$G$ is isomorphic to the direct limit of\/
$G_0\stackrel{\phi_0}{\rightarrow}G_1\stackrel{\phi_1}{\rightarrow}
G_2\stackrel{\phi_2}{\rightarrow}\dots$\textup.
\end{property}

\begin{proof}
Obviously, $G$ is isomorphic to the direct limit of
$G_0\stackrel{\phi_0}{\rightarrow}G_1\stackrel{\phi_1}{\rightarrow}
G_2\stackrel{\phi_2}{\rightarrow}\dots$.
The groups $G_0$, $G_1$, $G_2$, \dots are hyperbolic by
the Corollary \ref{corollary.proposition:1.(5.3)} of
Proposition~\ref{proposition:(5.3)}.
\end{proof}

It is easy to see that every recursively presented group with
a finite regular file basis has solvable word problem.
(Note that if a finitely generated group has a recursive presentation
with a countably infinite alphabet,
then it also has a recursive presentation with a finite alphabet.)
In particular, if the group $G$ is recursively presented,
then it has solvable word problem.
Generally, a recursively presented finitely generated group $H$ has
solvable word problem if an only if
there exists a recursive set of group words $X$ in a finite alphabet
$\mathfrak B$ such that
relative to some mapping of $\mathfrak B$ to $H$,
every element of $H$ has a unique representative in~$X$.
In the case of $G$, the set of all regular group words over
$\mathfrak A$ may
be taken as such a recursive set of representatives.

\begin{property}[conditional]
\label{property:5}
If\/ $G$ is recursively presented\textup,
then it has solvable conjugacy problem\textup.
\end{property}

\begin{proof}
Assume $G$ is recursively presented.
Then, as noted above, $G$ has solvable word problem.

Let $q$ be an integer such that $q\ge1/(1-2\mu)$.

Here is a description of an algorithm that solves the conjugacy
problem in~$G$:

\begin{quotation}
\raggedright
\tt
\begin{description}
\item[Input]
	group words $u$ and $v$ (in the alphabet~$\mathfrak A$).
\item[Step~1]
	If $[u]\ne1$ and $[v]\ne1$, go to the next step;
	otherwise, test whether $[u]=[v]$.
	If yes, output ``Yes'', if no, output ``No''.
	Stop.
\item[Step~2]
	Produce the (finite) set $\mathcal S$ of all
	group words $s$ over $\mathfrak A$ such that $[s]=1$ and
	$|s|\le q(|v|+|u|)$.
\item[Step~3]
	Produce a (finite) set $D$ of up to isomorphism all annular diagrams
	over $\langle\,\mathfrak A\,\|\,\mathcal S\,\rangle$
	with up to $q(|v|+|u|)$ edges.
\item[Step~4]
	If there is $\Delta\in D$ such that
	$u$ is the label of a representative
	of one of the contours of $\Delta$ and
	$v^{-1}$ is the label of some representative of
	the other, then output
	``Yes''; otherwise, output ``No''.
	Stop.
\end{description}
\end{quotation}

On an input $(u,v)$, this algorithm gives ``Yes'' as the output if
and only if $[u]$ and $[v]$ are conjugate in $G$; otherwise, it gives ``No.''
Below is an outline of a proof.

It is clear that the algorithm gives ``Yes'' only if $[u]$ and $[v]$
are conjugate in~$G$.
It is also clear that on every input it gives either ``Yes'' or ``No,''
and terminates.

Now, let $u$ and $v$ be arbitrary group words representing conjugate
elements of $G$, and consider $(u,v)$ as an input to the algorithm.

Case~1: $[u]=1$ or $[v]=1$.
Then $[u]=1=[v]$, and therefore the algorithm gives ``Yes'' on Step~1.

Case~2: $[u]\ne1$ and $[v]\ne1$.
Let $\Delta$ be a weakly reduced annular disc diagram over
$\langle\,\mathfrak A\,\|\,\mathcal R\,\rangle$ such that
$u$ and $v^{-1}$ are the labels of some representatives
of its two contours
%$u$ is the label of a representative
%of one of the contours of $\Delta$ and
%$v^{-1}$ is the label of a representative
%of the other contour of $\Delta$
(such $\Delta$ exists by Lemma V.5.2 in \cite{LyndonSchupp:2001:cgt}).
Let $\mathcal S$ be the set of group words produced by the algorithm on Step~2.
Every element of $\mathcal R$ of length at most $q(|v|+|u|)$ is in $\mathcal S$.
Let $D$ be the set of annular diagrams produced by the algorithm on Step~3.
By Proposition \ref{proposition:(5.1)} and the Main Lemma,
$$
|u|+|v|\ge(1-2\mu)\!\sum_{\Pi\in\Delta(2)}\!|\partial\Pi|
\ge\frac{1}{q}\sum_{\Pi\in\Delta(2)}\!|\partial\Pi|
$$
(recall that a special selection exists on every diagram
over $\langle\,\mathfrak A\,\|\,\mathcal R\,\rangle$).
Therefore, the length of the contour of every face of $\Delta$ is
at most $q(|u|+|v|)$, and the total number of edges of $\Delta$ is
\begin{multline*}
\|\Delta(1)\|
=\frac{1}{2}\Bigl(\sum_{\Pi\in\Delta(2)}\!|\partial\Pi|+|u|+|v|\Bigr)
\le\frac{1}{2}\bigl(q(|u|+|v|)+|u|+|v|\bigr)\\
=\frac{q+1}{2}(|u|+|v|)\le q(|u|+|v|).
\end{multline*}
Hence, $\Delta$ is isomorphic to some element of $D$,
and the algorithm returns ``Yes'' on Step~4.
\end{proof}

\begin{proof}[Proof of the theorem\/ \textup(see the Introduction\/\textup)]
Let $N$ be an integer such that $\lambda_1nN\ge1$.\break
Let the construction in Section~\ref{section:group_construction}
be carried out under
the following two additional conditions:
\begin{enumerate}
\item
	the order imposed on the set of all reduced group words is deg-lex;
\item
	for every $i$ for which $r_i$ is defined, $m_i=N|w_i|+i$.
\end{enumerate}
Evidently, these conditions are consistent with the rest.

By Proposition~\ref{proposition:(5.4)},
the group $G$ is recursively presented.
Properties \ref{property:1}--\ref{property:5} show that
up to isomorphism $G$ is a desired group.
($G$ is \emph{isomorphic\/} to a direct limit of a sequence
of hyperbolic groups with respect to a family of surjective homomorphisms,
but a desired group must \emph{be\/} such a limit.)
\end{proof}

%%%%%%%%%%%%%%%%%%%%%%%%%%%%%%%%%%%%%%%%%%%%%%%%%%%%%%%%%%%%%%%%%%%%%%%%%%%
 %<------------\|- Properties of the group
%% AMS-LaTeX
%% Paper 2 - Version GAMMA
%% Part 06
%% (Comments)
%%%%%%%%%%%%%%%%%%%%%%%%%%%%%%%%%%%%%%%%%%%%%%%%%%%%%%%%%%%%%%%%%%%%%%%%%%%

%!TEX root = paper2-2005.tex

\section{Comments}
\label{section:comments}

The group $G$ and its presentation 
$\langle\,\mathfrak A\,\|\,\mathcal R\,\rangle$
mentioned in this section are defined in 
Section~\ref{section:group_construction}.

\begin{proposition}
\label{proposition:(6.1)}
No hyperbolic group can be used as an example to demonstrate
the negative answer to Bludov's question\textup.
\end{proposition}

\begin{proof}
Suppose a hyperbolic group $H$ is an example demonstrating 
the negative answer.
Then $H$ is \emph{boundedly generated\/}
(is the product of a finite sequence of its cyclic subgroups)
and not virtually polycyclic.
In particular, $H$ is non-elementary.
Corollary~4.3 in \cite{Minasyan:2004:opqshg} states that every 
boundedly generated hyperbolic group
is elementary.
(This fact also follows from Corollary~2 
in~\cite{Olshanskii:1993:orhGshg}.)
This gives a contradiction.
\end{proof}

\begin{corollary.proposition}
\label{corollary.proposition:1.(6.1)}
The group\/ $G$ is not finitely presented\textup.
\end{corollary.proposition}

\begin{proof}
Suppose $G$ is finitely presented.
Then the presentation $\langle\,\mathfrak A\,\|\,\mathcal R\,\rangle$ 
is finite.
By Corollary \ref{corollary.proposition:1.(5.3)} 
of Proposition~\ref{proposition:(5.3)}, $G$ is hyperbolic.
This contradicts Proposition~\ref{proposition:(6.1)}.
\end{proof}

Thus, the set $\mathcal R$ and the sequences
$\{w_i\}_{i=1,\dots}$, $\{m_i\}_{i=1,\dots}$, $\{r_i\}_{i=1,\dots}$
constructed in Section~\ref{section:group_construction} are infinite.
\par

%%%%%%%%%%%%%%%%%%%%%%%%%%%%%%%%%%%%%%%%%%%%%%%%%%%%%%%%%%%%%%%%%%%%%%%%%%%
 %<------------\|- Comments
%% AMS-LaTeX
%% Paper 2 - Version GAMMA
%% Part 99
%% (Acknowledgements)
%%%%%%%%%%%%%%%%%%%%%%%%%%%%%%%%%%%%%%%%%%%%%%%%%%%%%%%%%%%%%%%%%%%%%%%%%%%

%!TEX root = paper2-2005.tex

\section*{Acknowledgements}

The author thanks Alexander Ol'shanskiy for suggesting the problem
and engaging in fruitful discussions during the course of this work.
\par

%%%%%%%%%%%%%%%%%%%%%%%%%%%%%%%%%%%%%%%%%%%%%%%%%%%%%%%%%%%%%%%%%%%%%%%%%%%
 %<------------\|- Acknowledgements
%\newpage
%% AMS-LaTeX
%% Paper 2 - Version GAMMA
%% Bibliography
%% (References)
%%%%%%%%%%%%%%%%%%%%%%%%%%%%%%%%%%%%%%%%%%%%%%%%%%%%%%%%%%%%%%%%%%%%%%%%%%%

%!TEX root = paper2-2005.tex

\nocite{Bludov:1995:nbg-rus,kn13:1995:ktnvtg-rus,Olshanskii:1989:gosg-rus}

\bibliographystyle{amsalpha}
\bibliography{\noexpand~/Documents/bib}

%%%%%%%%%%%%%%%%%%%%%%%%%%%%%%%%%%%%%%%%%%%%%%%%%%%%%%%%%%%%%%%%%%%%%%%%%%% %<------------\|- Bibliography

\end{document}